\documentclass[11pt,êeqno]{article}
\usepackage{amsfonts}
\usepackage[tbtags]{amsmath}
\topskip=\baselineskip%
\newcommand{\nty}{n \to \infty}

\newcommand{\lr}{\left(}
\newcommand{\lb}{\left[}
\newcommand{\rb}{\right]}
\newcommand{\rr}{\right)}
\newcommand{\lf}{\left\{}
\newcommand{\rf}{\right\}}
\newcommand{\lv}{\left|}
\newcommand{\rv}{\right|}
\newcommand{\lp}{\left.}
\newcommand{\rp}{\right.}
\newcommand{\kn}{k_n}
\newcommand{\mn}{m_n}
\newcommand{\eel}{\end{lemma}}

\newcommand{\xin}{X_{i:n}}

\newcommand{\uin}{U_{i:N}}

\newcommand{\xia}{\xi_{\alpha}}
\newcommand{\qa}{\xi_{\alpha_n}}
\newcommand{\qb}{\xi_{1-\beta_n}}
\newcommand{\qan}{\xi_{{\alpha_n n}:n}}

\newcommand{\an}{\alpha_n}
\newcommand{\bn}{\beta_n}

\newcommand{\nau}{N_{\alpha,u }}

\newcommand{\na}{N_{\alpha}}

\newcommand{\alp}{\alpha}

\newcommand{\aln}{\alpha_n}
\newcommand{\bln}{\beta_n}

\newcommand{\be}{\beta}

\newcommand{\si}{\sigma}

\numberwithin{equation}{section}

\newtheorem{theorem}{Theorem}[section]

\newtheorem{lemma}{Lemma}[section]

\newtheorem{corollary}{Corollary}[section]

\newtheorem{remark}{Remark}[section]

\title{{\huge \sf Second Order Approximations for Slightly Trimmed Sums}
}
\author{{\bf N.V.
Gribkova}\footnote{St.Petersburg~State~University,
Mathematics~and~Mechanics~Faculty, 198504, St.Petersburg,
Stary~Peterhof, Universitetsky~pr.~28,~Russia;
E-mail:~nv.gribkova@gmail.com}\ , \ {\bf R.
Helmers}\footnote{Centre~for~Mathematics~and~Computer~Science,
P.O.Box~94079,~1090~GB~Amsterdam,~The~Netherlands;
E-mail:~R.Helmers@cwi.nl}}
\date{
}
\begin{document}
\maketitle 

\begin{abstract}
{\small \it We investigate the second order asymptotic behavior of
trimmed sums $T_n=\frac 1n \sum_{i=\kn+1}^{n-\mn}\xin$, where
$\kn$, $\mn$ are sequences of integers, $0\le \kn < n-\mn \le n$,  such that
$\min(\kn, \mn) \to \infty$, as $\nty$, the $\xin$'s denote the order
statistics corresponding to a sample $X_1,\dots,X_n$ of $n$ i.i.d. random variables.
In particular, we focus on the case of
slightly trimmed sums with vanishing trimming percentages, i.e. we assume that  $\max(\kn,\mn)/n\to 0$, as $\nty$, and heavy tailed distribution $F$, i.e. the common distribution of the observations $F$ is supposed to have an~infinite variance.

We derive optimal bounds of Berry -- Esseen type of the order 
$O\bigl(r_n^{-1/2}\bigr)$, $r_n=\min(\kn,\mn)$,
for the normal
approximation to $T_n$ and, in addition, establish one-term  expansions of the
Edgeworth type for slightly trimmed sums  and
their studentized versions.

Our results supplement previous work on first order approximations
for slightly trimmed sums by Csorgo, Haeusler  \& Mason (1988) and
on second order approximations for (Studentized) trimmed sums with
fixed trimming percentages by Gribkova \& Helmers~(2006, 2007).}
\end{abstract}
\begin{quote}
\end{quote}

\section{Introduction and main results}
\label{imtro}
Let  $X_1,X_2,\dots $ be a sequence of independent identically
distributed (i.i.d.) real-valued nondegenerate random variables
(r.v.) with common distribution function ($df$) $F$, and for each
integer $n \ge 1$ let \ $X_{1:n}\le \dots \le X_{n:n}$ denote the
order statistics based on the sample $X_1,\dots ,X_n$. Introduce
the left-continuous inverse function $F^{-1}$ defined as
$F^{-1}(u)= \inf \{ x: F(x) \ge u \}$, \ $0<u\le 1$, \
$F^{-1}(0)=F^{-1}(0^+)$, and let $F_n$ and $F_n^{-1}$ denote the
empirical $df$ and its inverse respectively.

Define  the population truncated mean and variance functions
\begin{equation}
\label{1.1} \begin{split}
\mu (u,1-v)&=\int_u^{1-v} F^{-1}(s) \,d\, s  \ \ ,\\
\si^2 (u,1-v)&=\int_u^{1-v}\int_u^{1-v} (s\wedge t - st)\, d\,
F^{-1}(s)\, d\, F^{-1}(t),
\end{split}
\end{equation}
where $0\le u<1-v\le 1$, and $s\wedge t = \min (s,t)$. Note that $\si^2 (0,1)$ equals the variance of $X_1$ whenever $\boldsymbol{E}X_1^2$ is finite.

Let $\kn$ and $\mn$ be sequences of integers such that $0\le \kn <
n-\mn \le n$, and $\kn \wedge \mn \to \infty$, as $\nty$. Put
$\aln=\kn/n$, \ $\bln=\mn/n$.

Consider the trimmed sum given by
\begin{equation}
\label{1.2} T_n=\frac 1n \sum_{i=\kn+1}^{n-\mn}\xin
=\int_{\aln}^{1-\bln}F_n^{-1}(u)\, d\, u.
\end{equation}

The first order asymptotic properties of trimmed sums and slightly
trimmed sums (i.e. $\aln \vee\bln \to 0$) were investigated by many authors (cf.~\cite{s}, \cite{chm}, \cite{cm}, \cite{griffin}
and references therein). In particular in Cs\"orgo et al.~\cite{chm}
a~necessary and sufficient condition for the existence of $\{a_n\}$, $\{b_n\}$ such that
the distribution of the properly normalized slightly trimmed sum $a_n^{-1}(T_n-b_n)$ tends to the standard normal law was obtained, and (using a~different approach than in~\cite{chm}) Griffin and Pruitt~\cite{griffin} derived an~equivalent $iff$ condition for asymptotic normality of $T_n$.  In Griffin and Pruitt~\cite{griffin} the class of all subsequential limit laws for the sequences of slightly trimmed sums $a_n^{-1}(T_n-b_n)$ was characterized and sufficient conditions were given for $F$ to be in the domain of partial attraction of a given law from this class.
The members of this class are of the form $\tau N_1+f(N_2)-g(N_3)$, where $N_1$, $N_2$, $N_3$ are independent $N(0,1)$, $\tau \ge 0$ and $f$ and $g$ are arbitrary nondecreasing convex functions.
Both in \cite{chm} and in ~\cite{griffin} a~classical result by Stigler~\cite{s} for
the trimmed mean with fixed trimming percentages was extended to the case that
the fraction of trimming data is vanishing when $n$ gets large.

The second order asymptotic properties  of trimmed sums with fixed trimmed percentages
 were investigated  by Gribkova and
Helmers~\cite{gh2006}-\cite{gh2007}: the validity of the one-term Edgeworth
expansion (EE) for a (Studentized) trimmed mean and bootstrapped
trimmed mean were established and simple explicit formulas of the
first leading terms of these expansions were found.
(We note in passing that in Helmers et
al.~\cite{hjq} a saddlepoint approximation -- a completely different
way of approximating $df$ of the trimmed mean with fixed trimming percentages
accurately -- was obtained.)

Here we extend the result in \cite{gh2006} and establish second order asymptotic properties to slightly trimmed means and their Studentized versions. Except for a~non optimal Berry -- Esseen type bound for slightly trimmed means in \cite{en}, to the best of our knowledge, no such second order results are available in the literature.
In this article we focus on the case of heavy-tailed distributions, i.e. we assume that  $\si^2(0,1)=\infty$ (cf.~\eqref{1.1}). We refer to Remark~1.1 for a~detailed discussion of the different behavior of our second order approximations in the case of a~heavy tailed respectively light tailed distribution~$F$.

To begin with we shall obtain bounds of Berry -- Esseen type for
the normal approximation to  $T_n$ under a~weak condition on the density, assuming its existence in the tails
of the distribution of the observations, we also show that the bounds we give in this paper, namely
$O((\kn \wedge \mn)^{-1/2})$, in absence of any moment assumptions  are of the best possible order.

Secondly we will supplement our results on the rates of convergence towards normality by deriving the one-term expansions of the Edgeworth type for slightly trimmed
 sums and Studentized slightly
trimmed sums and obtain simple explicit formulas for these
expansions. In a way we in particular refine the first order limit results of Cs\"orgo et al.~\cite{chm} by establishing more accurate second order approximation of Edgeworth type for slightly trimmed means with vanishing trimming percentages.

We show  that  the first leading term of our one-term expansion (in absence of symmetry) has the exact order $(\kn \wedge \mn)^{-1/2}$  (cf.~\eqref{opt}, Remark~1.1), when the density is regular varying in the tails with index $\rho=-(1+\gamma), \  0<\gamma<2$, (cf.~Bingham et.al~\cite{bingham} on the topic regular variation and Borovkov and Mogulskii~\cite{bormogu} for assumptions similar to our condition $[R]$ on~p.~\pageref{r1}), which directly imply optimality of our Berry -- Esseen  type bounds.  (When the underlying distribution has finite second moment, the order of the bound can be improved to $O(n^{-3(1/2-1/\gamma)})$ when  $2<\gamma<3$, and to $O(n^{-1/2})$ when $3\le \gamma$; $-(1+\gamma)$ is the index of regularity of the density in the tails (cf.~condition $[R]$ on~p.~\pageref{r1}). We will pursue this topic  elsewhere).

Similarly as in~\cite{gh2006}-\cite{gh2007}, our method of
proof is based on a~stochastic approximation of a~(slightly) trimmed
sum by a~$U$-statistic of degree two with a~kernel depending on~$n$.
We use also a Bahadur--Kiefer type approximation (cf.~section~\ref{Bahadur} ).

We conclude this introduction by noting that the case of heavy-tailed distribution  we focus on in this article is interesting in particular due to the following statistical motivation: suppose that $\boldsymbol{E}|X_1|<\infty$ and that we are interested in estimating of $\boldsymbol{E}X_1$, whereas
 $\si^2(0,1)=\infty$. The trimmed mean with fixed trimming percentages (robust estimate) is not consistent in absence
of symmetry of the underlying distribution, but the slightly trimmed mean $T_n$ tends to $\boldsymbol{E}X_1$ a.s., so it is a~consistent estimator. The next issue one may want to consider is interval estimation. Fortunately, the suitably normalized (or studentized)   $T_n$ has a~standard normal asymptotic distribution. The rate of convergence in case of the heavy tailed $F$ can be rather slow (cf. Theorem~1.2 and Corollary 1.1 below). However, if we know a second order asymptotic approximation to a~$df$ of the normalized $T_n$ (cf. Theorem~1.5) and of the studentized $T_n$ (cf.~Theorem~1.7)  correcting the bias and skewness, we can improve the standard normal approximation to an~approximation having smaller remainder.

The paper is organized as follows: in Section~1, we formulate sets
of conditions and state our main results on Berry -- Esseen type
bounds and the Edgeworth type expansions for a~normalized (slightly) trimmed sums and for its
Studentized versions. In Section~2, we state and prove
the~auxiliary results on the $U$-statistic approximation for $T_n$ and
for the plug-in estimate of its asymptotic variance. The proofs of the
main results are relegated to Section~3. In Section~\ref{Bahadur}, we state and
prove two Bahadur--Kiefer type lemmas,  which we use in our proofs,,  in particular, lemma~\ref{lem_A.2} provides a representation for a sum of order statistics lying between  the $\an$-th
population quantile and the corresponding empirical quantile. A lemma used  in the proofs of the Bahadur--Kiefer type results is relegated to the Appendix.

Define the $\nu$-th quantile of $F$ by $\xi_{\nu}=F^{-1}(\nu)$,
$0<\nu<1$, and let $W_i(n)$, \ $i=1,\dots,n$, denote $X_i$
Winsorized outside of $(\qa,\qb]$, that is
\begin{equation}
\label{1.3} W_i(n)=\qa \vee ( X_i \wedge \qb),
\end{equation}
where $s\vee t=\max(s,t)$. Define the quantile function
\begin{equation}
\label{1.4} Q_n(u)=\qa \vee (F^{-1}(u) \wedge \qb),
\end{equation}
and the first two cumulants of $W_i(n)$:
\begin{equation}
\label{1.5} \mu_{W_{(n)}}=\int_0^1Q_n(u) d\, u \ , \quad
\si^2_{W_{(n)}}=\int_0^1(Q_n(u) -\mu_{W_{(n)}})^2d\, u \ .
\end{equation}
Note that $\si^2_{W_{(n)}}=\si^2(\an,1-\bn)$ (cf. \eqref{1.1}),
and its square root is a suitable scale parameter for $T_n$ when
establishing its asymptotic normality (cf. Cs\"org\H{o} et al.~\cite{chm},
see also Griffin and Pruitt~\cite{griffin}). We will suppose throughout
this article that $\liminf_{\nty}\si_{W_{(n)}}>0$ (i.e. $\qa
\neq \qb$ for all sufficiently large $n$).

Define four numbers
\begin{equation}
\label{seq}
\begin{split}
a_1&= \liminf_{\nty}\, \an , \qquad  \qquad a_2= \limsup_{\nty}\, \an , \\
b_1&= \liminf_{\nty}\, (1-\bn), \qquad  b_2= \limsup_{\nty}\,
(1-\bn) ,
\end{split}
\end{equation}
where $0\le a_1\le a_2$, \ $b_1\le b_2\le 1$, and suppose that
$a_2<b_1$.

We will assume throughout this article that the following
smoothness condition is satisfied.

\medskip \noindent
$[A_1].$ \ \ {\it There exist two open sets \ $U_a$, $U_b\subset
(0,1)$ such that $F^{-1}$ is differentiable in $U=U_a \cup U_b$, and
\begin{equation}
\label{1.6}
\begin{array}{cccccccccc}
 & (0,\varepsilon),\ \ \ &\text{if }\ \ 0=a_1=a_2,&&(1-\varepsilon,1),& \text{if }\ \ b_1=b_2=1,\\
U_a \ \supset & (0,a_2], \ \  &\text{if }\ \ 0=a_1<a_2,&\qquad U_b \ \supset &[b_1,1),\ \ \ \ &\text{if }\ \  b_1<b_2=1,\\
& [a_1,a_2], \ &\text{if }\ \   0<a_1\le a_2,& &[b_1,b_2],\ \ \  & \text{if }\ \ b_1\le b_2<1,\\
\end{array}
\end{equation}
(with some $0<\varepsilon  <1$ in cases given in the first lines of
\eqref{1.6}), i.e. the density $f=F'$ exists and is positive in  $F^{-1}(U)$.}

Define two sequences:
\begin{equation}
\label{qu} q_{\an}=\frac 1{\sqrt{n}\, \si_{W_{(n)}}} \, \frac
{\an}{f(\qa)},\qquad q_{\bn}=\frac 1{\sqrt{n}\, \si_{W_{(n)}}} \,
\frac {\bn}{f(\qb)}.
\end{equation}
We note that it is a simple consequence of $[A_1]$ that these quantities are well defined for
all sufficiently large $n$. The same remark also applies to some
other quantities we introduce below.

Our second assumption is:

\medskip\noindent$[A_2].$
\qquad \qquad \qquad $q_{\an}\ \vee \ \ q_{\bn} \longrightarrow
0$, \quad {\it as} \quad $\nty$.

\bigskip
Note that $[A_2]$ holds true if $\si^2=\si^2(0,1)<\infty$ and
$\frac{\an}{f(\qa)}\vee \frac{\bn}{f(\qb)}=o(\sqrt{n})$. Moreover,
$[A_2]$ is also satisfied if $\si^2=\infty$  and
$\frac{\sqrt{\an}}{|\qa| f(\qa)}\ \ \vee \frac{\sqrt{\bn}}{|\qb|
f(\qb)}=o(\sqrt{n})$, because, due to Lemma~2.1 of Cs\"org\H{o} et
al.~\cite{chm}, for any quantile function $F^{-1}$:
\begin{equation}
\label{1.7} \limsup_{u,v\downarrow 0}\frac{u\, (F^{-1}(u))^2+v\,
(F^{-1}(1-v))^2}{\si^2(u,1-v)}<\infty\, .
\end{equation}

In a way relation \eqref{1.7}  will be  crucial for our purposes.
Note first of all that in the special case that the second moment
of $F$ is assumed to be finite (cf.~Theorem~1.3) the upper limit 
in \eqref{1.7} is not only bounded but is in fact equal
to zero. This simple fact is at the basis of the slightly better
rates obtained in Theorem~1.3 in comparison with the rate
established in the more general Theorem~1.2.

Let $h$ be a real-valued function defined on the set $F^{-1}(U)$
(cf. \eqref{1.6}). Take an arbitrary $0<B<\infty$ and for all sufficiently large $n$ define
\begin{equation}
\label{1.8}
\begin{split}
\Psi_{\an,h}(B)&=\sup_{|t|\le B}\lv h\circ F^{-1}\Bigl( \an
+t\sqrt{\frac{\an \ln \kn}{n}}\Bigr) -h\circ F^{-1}\Bigl(\an\Bigr)\rv,\\
\Psi_{1-\bn,h}(B)&=\sup_{|t|\le B}\lv h\circ F^{-1}\Bigl( 1- \bn
+t\sqrt{\frac{\bn \ln \mn}{n}}\Bigr) -h\circ
F^{-1}\Bigl(1-\bn\Bigr)\rv,
\end{split}
\end{equation}
where $h\circ F^{-1}(u)=h\lr F^{-1}(u)\rr$. Note that $\an
+t\sqrt{\frac{\an \ln \kn}{n}}=\an \lr 1 +t\sqrt{\frac{\ln
\kn}{\kn}}\rr =\an (1+o(1))$, and $1-\bn +t\sqrt{\frac{\bn \ln
\mn}{n}}=1-\bn (1+o(1))$, as $\nty$. In particular, this implies that the two functions introduced in~\eqref{1.8} are
well-defined for all sufficiently large~$n$.

We will use in what follows the auxiliary functions:
$\Psi_{\nu_n,x}(B)$, $\Psi_{\nu_n,\frac 1{f(x)}}(B)$
$\Psi_{\nu_n,\frac x{f(x)}}(B)$, $\nu_n=\an,\, 1-\bn$,
corresponding to $h(x)=x$, $1/f(x)$ and $x/f(x)$ in~\eqref{1.8}. It is easy to
see in any case that the following inequalities are valid:
\begin{equation}
\label{eq_psi}
\begin{split}
\Psi_{\an,\frac x{f(x)}}(B)&\le \an
B\Bigl(\frac{\ln\kn}{\kn}\Bigr)^{1/2}\lr \frac 1{f(\qa)}+
\Psi_{\an,\frac 1{f(x)}}(B)\rr^2  + |\qa|\Psi_{\an,\frac 1{f(x)}}(B),\\
\Psi_{1-\bn,\frac x{f(x)}}(B)&\le \bn
B\Bigl(\frac{\ln\mn}{\mn}\Bigr)^{1/2}\lr \frac 1{f(\qb)}+
\Psi_{1-\bn,\frac 1{f(x)}}(B)\rr^2 +  |\qb|\Psi_{1-\bn,\frac 1{f(x)}}(B).
\end{split}
\end{equation}
These inequalities will be especially useful in the proof of
Lemma~2.2 in Section~2.

 Our third assumption is

\medskip \noindent
$[A_3]$ \ {\it For every $0<B<\infty$
\begin{equation*}
\frac {\an}{\sqrt{n}\, \si_{W_{(n)}}} \, \Psi_{\an, \frac
1{f(x)}}(B) \longrightarrow 0, \qquad \  \frac {\bn}{\sqrt{n}\,
\si_{W_{(n)}}} \, \Psi_{1-\bn, \frac 1{f(x)}}(B)\longrightarrow 0,
\end{equation*}
as $\nty$}.

Define the $df$ of the normalized $T_n$ by
\begin{equation}
\label{1.9} F_{T_n}(x)=\boldsymbol{P}\lr \si_{W_{(n)}}^{-1}n^{1/2}\lr \frac
{}{} T_n-\mu(\an,1-\bn)\rr \le x \rr .
\end{equation}

First we show that the conditions $[A_1]$ -- $[A_3]$ together
yields
\begin{equation}
\label{clt}\sup_{x\in R}|F_{T_n}(x)-\Phi (x)|=o(1)\, , \ \
\text{as} \ \ \nty \, ,
\end{equation}
where $\Phi$ is standard normal $df$. To check this we verify that
the $iff$ conditions of asymptotic normality of the trimmed sum
$T_n$ (cf. Cs\"org\H{o} et al.~\cite{chm},Theorem~4, p.~677) are automatically satisfied
whenever our conditions $[A_1]$ - $[A_3]$ hold true. Consider the
first auxiliary function defined on page 674 of Cs\"org\H{o} et
al.~\cite{chm}, which corresponds to the trimming of the $k_n$
smallest observations on our sample of size n ; the treatment of
the second auxiliary function on p. 674 of the same paper, which deals with the trimming of the $\mn$ largest observations, is
similar and therefore omitted. For \eqref{clt} to hold we must
verify that for every $c\in \mathbb{R}$
\begin{equation}
\label{1.10} Q_{1,n}(c)\to 0, \ \ \nty,
\end{equation}
where
\begin{equation*}
Q_{1,n}(c)= \begin{cases} \frac{(\an)^{1/2}}{\si_{W_{(n)}}}\lf
F^{-1}\Bigl( \an +c\sqrt{\frac{\an}{n}}\Bigr) -F^{-1}\Bigl(\an\Bigr)\rf\ , \\
\qquad \qquad \qquad \qquad \qquad \qquad \ \qquad |c|\le \frac 12 \sqrt{\an n}\ , \\
Q_{1,n}\bigl(- \frac 12 \sqrt{\an n}\bigr), \ \ \ \ \qquad -\infty < c <-\frac 12 \sqrt{\an n}\ ,\\
Q_{1,n}\bigl(\frac 12 \sqrt{\an n}\bigr), \ \ \ \ \  \frac 12
\sqrt{\an n}<c<\infty \ ,
\end{cases}
\end{equation*}
(cf. Cs\"org\H{o} et al.~\cite{chm}). Note that $\an n=\kn \to\infty$,
and for each $c\in \mathbb{R}$ and all sufficiently large $n$ we
have  $|c|<\frac 12 \sqrt{\an n}$ , and $\an
+c\sqrt{\frac{\an}{n}}=\an(1+c\sqrt{\frac{1}{\kn}})$ belongs to
 the set $U_a$ (cf. \eqref{1.6}). So we have
\begin{equation}
\label{1.13}
 Q_{1,n}(c)=\frac{(\an)^{1/2}}{\si_{W_{(n)}}}\, c\,
\sqrt{\frac{\an}{n}}\frac 1{f\Bigl( F^{-1}\Bigl( \an +\theta
c\sqrt{\frac{\an}{n}}\Bigr) \Bigr)}
\end{equation}
for some $0<\theta<1$, and the quantity \eqref{1.13} in absolute
value is less than
\begin{equation*}
 |c|\, \frac {\an}{\sqrt{n}\, \si_{W_{(n)}}} \Big(\frac 1{f(\qa)}
+\Psi_{\an,1/f(x)}(\theta |c|) \Big)\ ,
\end{equation*}
which tends to zero by $[A_2]$ and $[A_3]$, and \eqref{1.10}
follows.

Our conditions $[A_1]$ -- $[A_3]$ are slightly stronger than $iff$
conditions of asymptotic normality of Cs\"org\H{o} et al.~\cite{chm},
but these conditions enable us to establish a~bound for the error in the normal
approximation for the $df$ of $T_n$ .

We note in passing that Peng~\cite{peng} has shown that it is
impossible in general to replace the truncated mean $\mu(\an,
1-\bn)$ employed in \eqref{1.9} by the  ordinary mean of the
trimmed sum $\boldsymbol{E}T_n$ (which is always finite of course when $F$ and $1-F$ are regular varying
at minus and plus infinity respectively); centering by a truncated mean is
really needed to obtain a standard normal limit in \eqref{1.9}.

Here is our general result on the rate of convergence of the
distribution of a~properly normalized trimmed sum $T_n$ to the standard
normal law.
\begin{theorem}
Assume that the conditions $[A_1]$ and $[A_2]$ are satisfied. Then
\begin{equation}
\label{BE_1} \sup_{x\in \mathbb{R}}\lv F_{T_n}(x) - \Phi(x) \rv
\le \frac
A{\sqrt{n}}\bigl(\delta_{1,n}+\delta_{2,n}+\delta_{3,n}+\delta_{4,n}
\bigr)+ C\bigl(\kn^{-c}+\mn^{-c} \bigr),
\end{equation}
\begin{equation*}
\begin{split}
\delta_{1,n}&=\frac{\boldsymbol{E}\lv W_1(n)\rv^3}{\si^3_{W_{(n)}}}\,
,\qquad \ \ \delta_{2,n}=\frac 1{\si_{W_{(n)}}}\lr
\frac{\an}{f(\qa)}+\frac{\bn}{f(\qb)}\rr,\\
\delta_{3,n}&=\big(\an\big)^{1/3}\lr \frac {\an}{f(\qa)
\si_{W_{(n)}}}\rr^{5/3}+\big(\bn\big)^{1/3}\lr \frac {\bn}{f(\qb)
\si_{W_{(n)}}}\rr^{5/3},\\
\delta_{4,n}&=\frac 1{\si_{W_{(n)}}}\lr\an \ln\kn\,
\Psi_{\an,\frac 1{f(x)}}(B)+ \bn \ln\mn\, \Psi_{1-\bn,\frac
1{f(x)}}(B)\rr,
\end{split}
\end{equation*}
for every $c>0$, where $A,B,C>0$ are some constants, depending only
on $c$.
\end{theorem}

Note that  at the r.h.s. of \eqref{BE_1} we have:
$\frac1{\sqrt{n}}\delta_{1,n}=O(\frac 1{\sqrt{\kn}}+\frac
1{\sqrt{\mn}})$  in view of
\eqref{1.7} (cf.~Proof of Theorem~1.2, Section~3), $\frac1{\sqrt{n}}\delta_{2,n}=o(1)$ in view of $[A_2]$, $\frac1{\sqrt{n}}\delta_{3,n}=o(1)$ if we additionally assume that $q_{\an}=o(\kn^{-1/5})$,
 $q_{\bn}=o(\mn^{-1/5})$ and $\frac1{\sqrt{n}}\delta_{4,n}=o(1)$ if $\frac{\an}{\sqrt{n}\si_{W_{(n)}}}\Psi_{\an,\frac 1{f(x)}}(B)=o(\frac 1{\ln \kn})$, \ $\frac{\bn}{\sqrt{n}\si_{W_{(n)}}}\Psi_{1-\bn,\frac 1{f(x)}}(B)=o(\frac 1{\ln \mn})$.

 \medskip
 \noindent{\bf Example 1.1}. \ Let us consider an example, where the underlying distribution $F$ has super-heavy tails.  Let $F$ is such that $F(x)=1-F(x)=\frac 12\frac 1{(\ln |x|)^{\rho}}$, $\rho>0$, for all $x$: $|x|>x_0>0$ (cf.~\cite{griffin}). Simple computations on the quantities $\delta_{i,n}$ (it turns out that the term, corresponding $i=3$ is the largest one in this case) show that at the r.h.s. of \eqref{BE_1} we have a~bound of the order $O(\kn^{-d}+\mn^{-d})$ with some $0<d\le 1/2$ when
 $\limsup_{n\to\infty}(\kn^{-1}+\mn^{-1} )n^{\frac{5}{5+3\rho(1/2-d)}}<\infty$, and the bound of the order  $O(\kn^{-1/2}+\mn^{-1/2})$ is possible if and only if $\kn\asymp n$ and $\mn\asymp n$. We can obtain a~bound of the order $O(\kn^{-1/3}+\mn^{-1/3})$ (say) if we take $\kn$, $\mn$ $\sim \, n^{\frac{10}{10+\rho}}$

\medskip

To obtain more explicit bounds than the bound given in \eqref{BE_1} we need some more restrictive conditions. The
following assumption is somewhat stronger than $[A_2]$:

\medskip
\noindent$[A_2'].$ {\it \ Suppose that}
\begin{equation*}
\ \ \
\limsup_{\nty}\frac{\an^{3/2}}{\si_{W_{(n)}}f(\qa)}<\infty,\qquad
\limsup_{\nty}\frac{\bn^{3/2}}{\si_{W_{(n)}}f(\qb)}<\infty.
\end{equation*}
\noindent The latter condition implies
\begin{equation}
\label{qaqb}  q_{\an}=O\lr\frac{1}{\sqrt{\kn}}\rr , \qquad \ \
q_{\bn}=O\lr\frac{1}{\sqrt{\mn}}\rr, \ \ \text{\it as} \ \ \nty \,
.
\end{equation}

Note that in view of \eqref{1.7} condition $[A_2']$ holds true if
the following slightly  stronger condition is satisfied:

\medskip
\noindent$[A_2''].$ {\it Suppose that }
\begin{equation*}
\ \ \ \limsup_{\nty}\frac{\an}{|\qa|f(\qa)}<\infty,\qquad
\limsup_{\nty}\frac{\bn}{|\qb|f(\qb)}<\infty.
\end{equation*}

Note that in the case of a~slightly trimmed sum (i.e. when
$\an\vee\bn\to 0$) condition $[A_2'']$ is certainly satisfied when
$\limsup_{x\to -\infty}\frac{F(x)}{|x|f(x)}<\infty$ and
$\limsup_{x\to +\infty}\frac{1-F(x)}{xf(x)}<\infty$, and that the
latter requirement is true when the $df$ $F$ has a density for all
sufficiently large $|x|$, and $f=F'$ is regularly varying at the
infinity with index $\rho<-1$  (cf.~condition $[R]$, Corollary~1.1,~Theorems~1.5,~1.7).

The following condition is stronger than smoothness condition
$[A_3]$:

\medskip
\noindent$[A_3'].$ \ \ {\it For every $B>0$
\begin{equation*}
\ \ \ \Psi_{\an,\frac 1{f(x)}}(B) =O\lr \frac 1{f(\qa)\ln\kn}\rr\, ,
\qquad \Psi_{1-\bn,\frac 1{f(x)}}(B)=O\lr \frac 1{f(\qb)\ln\mn}\rr\,
,
\end{equation*}
}

Now we are in a position to state our second result of
Berry--Esseen type, which yields an~explicit upper bound of a~much
simpler form:
\begin{theorem}
Suppose that conditions $[A_1]$, $[A_2']$ and $[A_3']$ hold. Then
\begin{equation}
\label{BE_2} \sup_{x\in \mathbb{R}}\lv F_{T_n}(x) - \Phi(x) \rv
\le C\Bigl(\frac{1}{\sqrt{\kn}}\ +\ \frac{1}{\sqrt{\mn}}\Bigr),
\end{equation}
where $C$ is a positive constant not depending on \ $n$.
\end{theorem}

This result can be compared with an earlier result by Egorov \&
Nevzorov~\cite{en}, where a~non optimal bound of the order
$O\bigl(\frac{\ln\kn}{\sqrt{\kn}}\ +\
\frac{\ln\mn}{\sqrt{\mn}}\bigr)$ under stronger
conditions was obtained. In contrast, our bound \eqref{BE_2} is
sharp and yields an~optimal order bound of Berry--Esseen type for
slightly trimmed means when $F$ is, for instance, the Cauchy
distribution. The optimality of the bound in \eqref{BE_2} follows directly from our results on
the Edgeworth  type expansions and computations given in Remark~1.1 (cf.~\eqref{opt}).

Our next assertion concerns the case of a~slightly trimmed
mean for the special case when $\boldsymbol{E}X^2_1<\infty$, i.e.~the case of a~light tailed distribution.
\begin{theorem}
Suppose that $\si^2=\si^2(0,1)<\infty$, \ $\an\vee\bn\to 0$ as
$\nty$, and that the conditions $[A_1]$, $[A_2'']$ hold true. In
addition assume that for every $B>0$: \
\begin{equation}
\label{cond_psi}
\Psi_{\an,\frac
1{f(x)}}(B)=o\bigl((f(\qa)\ln\kn)^{-1}\bigr), \
\Psi_{1-\bn,\frac 1{f(x)}}(B)=o\bigl((f(\qb)\ln\mn)^{-1}\bigr),
\end{equation}
as $\nty$. \ Then
\begin{equation}
\label{BE_3} \sup_{x\in \mathbb{R}}\lv F_{T_n}(x) - \Phi(x) \rv
=o\Bigl(\frac{1}{\sqrt{\kn}}\ +\ \frac{1}{\sqrt{\mn}}\Bigr),
\end{equation}
as $\nty$.
\end{theorem}

This result --- i.e. the order bound \eqref{BE_3} --- applies for
instance to a $df$~$F$ with a~regular varying density $f$ which
behaves like $|x|^{-(3+\varepsilon)}$ (with some $\varepsilon >0$)
in the tails, so that the variance of $F$ is indeed finite. If we
take in addition -- by way of an example -- $k_n=m_n=[n^{1/2}]$
then we obtain a~sharper bound of order $o(n^{-1/4})$ instead of
$O(n^{-1/4})$ which would follow from the previous Theorem~1.2.
Moreover, if moments of higher order than $2$ are assumed to be
finite, it appears possible to establish the exact order of the
normal approximation error in \eqref{BE_3}, rather than asserting
only that the order of magnitude of the normal error in
\eqref{BE_3} is smaller than the one in \eqref{BE_2}. A~detailed
study of these exact rates, however, is outside the scope of the
present paper. The authors hope to pursue this matter elsewhere.

Next we obtain some consequences of Theorems~1.1,~1.2 and~1.3. Our
first corollary concerns the case of slightly trimmed sum when the
$df$ $F$ belongs to a domain of attraction of a stable law. Let
$RV_{\rho}^{\infty}$ be a class of regularly varying in the
infinity functions: $g\in RV_{\rho}^{\infty}$ $\Leftrightarrow$
$g(x)=|x|^{\rho}\, L(x)$, for $|x|>x_0$, with some $x_0>0$,
$\rho\in\mathbb{R}$, and $L(x)$ is a positive slowly varying function at infinity.
We will need the following regularity condition on the tails for the density $f$:

\medskip
\noindent$[R].$ \ {\it Suppose that $f\in RV_{\rho}^{\infty}$, \
where $\rho=-(1+\gamma)$, \ $\gamma >0$, \ and assume that
\begin{equation}
\label{r1} \lv f(x+\triangle x)-f(x)\rv=O\Bigl( f(x)\, \Big|\,
\frac{\triangle x}{x}\, \Big|\Bigr),
\end{equation}
when $\triangle x=o(|x|)$, as $|x|\to\infty$}.

\medskip
Note that \eqref{r1} holds true for $f$ if $\Big|\frac
{L(x+\triangle x)}{L(x)}-1\Big|=O\Bigl( \Big|\, \frac{\triangle
x}{x}\, \Big|\Bigr)\Bigr)$, as $|x|\to +\infty$, where $L$ is the
corresponding slowly varying function, and it is satisfied if $L$
is continuously differentiable for sufficiently large $|x|$ and
$|L'(x)|=O\lr\frac {L(x)}{|x|}\rr$, as $|x|\to +\infty$, which is
valid for instance when $L$ is some power of the logarithm. We
refer to Borovkov and Mogulskii~\cite{bormogu}, p.~568 for some conditions
closely related to ours.

The following corollary holds true for a slightly trimmed mean in
case of a~regular varying density:
\begin{corollary}
Suppose that $\an\vee\bn\to 0$, as $\nty$, condition $[A_1]$ holds
true for some $\varepsilon >0$, and the density $f$ satisfies $[R]$
with $0<\gamma \le 2$ on the set \ $F^{-1}(U)$. Then: \\
$(i)$ the bound
\eqref{BE_2} is valid; \\
$(ii)$ in addition if $\gamma=2$ and
$\si^2<\infty$ then also the sharper order bound \eqref{BE_3} holds true.
\end{corollary}

It is clear from our proofs (cf~Section~3) that the latter
assertion is valid as well if the density $f$ has different
indices of regularity near $-\infty$ respectively to
$+\infty$ (in particular, at least one of them can be greater than $2$), we
keep these two indices equal to each other for simplicity. Moreover this situation
corresponds to the important special case when $F$ belongs to a domain of
attraction of a~stable law.

Our second  corollary concerns the classical case when trimming
occur on the levels of the central order statistics. Let $a_1$,
$b_2$,  $U_a$ and $U_b$ are as in \eqref{seq}-\eqref{1.6}.
\begin{corollary}
Suppose that $0<a_1<b_2<1$, and  assume that the condition $[A_1]$
is satisfied. In the addition suppose that the density $f$
satisfies a~H\"{o}lder condition of degree~$d$ (for some $d>0$) on the sets
$F^{-1}(U_c)$, $c=a,\, b$. Then
\begin{equation}
\label{cor2} \sup_{x\in \mathbb{R}}\lv F_{T_n}(x) - \Phi(x) \rv
\le \frac C{\sqrt{n}}\, ,
\end{equation}
where $C>0$ is a constant, not depending on $n$.
\end{corollary}

Note that the smoothness assumptions imposed in Corollary 1.2 are
especially well suited  for obtaining our results on the Edgeworth type
expansions, which we will state and prove below. So the smoothness
assumption in corollary~1.2 is slightly excessive for obtaining of
the Berry -- Esseen type bound \eqref{cor2} (cf.~for instance,~\cite{gri}, where the optimal bound
was obtained under a somewhat weaker
smoothness assumption that $F^{-1}$ satisfies a Lipschitz
condition on the sets $U_a$ and  $U_b$, by an~application of Theorem~1.1 of van~Zwet~\cite{zwet} for
symmetric statistics).

\medskip
Next we will go one step further and establish one-term Edgeworth type expansions for $df$ of a~normalized and of a~Studentized slightly trimmed sum.

Define
\begin{equation*}
\gamma_{3,W_{(n)}}=\int_0^1(Q_n(u)-\mu_{W_{(n)}})^3\, d\,u,
\end{equation*}
where $Q_n(u)$, $\mu_{W_{(n)}}$ as in \eqref{1.4}--\eqref{1.5},
and put
\begin{equation}
\label{delta_2} \delta_{2,W_{(n)}}=-\alp_n^2\frac{\bigl(
\mu_{W_{(n)}} -\qa\bigr)^2} {f(\qa)}+\be_n^2\frac{\bigl(
\mu_{W_{(n)}} -\qb\bigr)^2} {f(\qb)}\, .
\end{equation}
Define two sequences of the real numbers
\begin{equation}
\label{la_12}
\lambda_{1_{(n)}}=\frac{\gamma_{3,W_{(n)}}}{\si^3_{W_{(n)}}}\
,\qquad
\lambda_{2_{(n)}}=\frac{\delta_{2,W_{(n)}}}{\si^3_{W_{(n)}}}\
\end{equation}

We establish the validity of the Edgeworth type expansion for the $df$
$F_{T_n}$ under conditions $[A_1]$-$[A_3]$. This expansion is
given by
\begin{equation}
\label{G_n} G_n(x)=\Phi(x)-\frac{\phi(x)}{6\sqrt{n}}\Bigl(\bigl(
\lambda_{1_{(n)}}+3\lambda_{2_{(n)}}\bigr)(x^2-1)+6\sqrt{n}\frac{b_n}{\si_{W_{(n)}}}\Bigr),
\end{equation}
where $\phi\,=\Phi'$, \ and \ \ $b_n=\frac 1{2\sqrt{n}}\lr
-\frac{\an(1-\an)}{f(\qa)}+\frac{\bn(1-\bn)}{f(\qb)}\rr$, $b_n$ is a~bias term which is present in the expansion
despite of the absence of any moment assumptions (cf.~\cite{gh2006}).

Note that if  $\an=\bn$ and the underlying distribution is symmetric, we have $G_n(x)\equiv \Phi(x)$ because the second term of the expansion is equal to zero in this case.

Similarly as when proving of Theorem 1.2 it is easy to check that if conditions $[A_1]$
 and $[A'_2]$ are satisfied the second term of $G_n(x)$ at the r.h.s. of \eqref{G_n} (for each fixed $x$) is a~magnitude of the order  $O\Bigl(\frac{1}{\sqrt{\kn}}\ +\ \frac{1}{\sqrt{\mn}}\Bigr)$ as $\nty$.
 And under some proper conditions (cf.~Theorem~1.5 and Corollary~1.3) the remainder in approximating of the $df$ of the normalized slightly trimmed sum by its expansion $G_n(x)$ is of the Bahadur type order $O\lr
\frac{(\ln\kn)^{5/4}}{\kn^{3/4}}+\frac{(\ln\mn)^{5/4}}{\mn^{3/4}}
\rr$, as $\nty\, $.

 Some simple computations show that in case of underlying distribution  $F$  considered in Example~1.1 ($F$ has no finite moments) $|G_n(x)-\Phi(x)|\asymp\frac 1{\sqrt{\kn}}(\frac n{\kn})^{\frac 1{\rho}}+\frac 1{\sqrt{\mn}}(\frac n{\mn})^{\frac 1{\rho}}$, $x\in R$. So, $|G_n(x)-\Phi(x)|\asymp \frac{1}{\sqrt{\kn}}\ +\ \frac{1}{\sqrt{\mn}}$ if and only if $\kn \asymp n$ and $\mn \asymp n$.

\medskip
 \noindent{\bf Remark~1.1} \ \ Let us investigate the order of magnitude of the various terms appearing in
 the Edgeworth type correction \eqref{G_n}. Two of these terms are correcting for skewness, let us denote them by
 $t_{j,n}=\frac 1{\sqrt{n}}\lambda_{j_{(n)}}$, $j=1,2$, respectively, while a third term is correcting for the bias present, which we denote by $t_{3,n}=b_n/\si_{W_{(n)}}$. Suppose now that $\max(\an,\bn)\to 0$, condition $[A_1]$ is satisfied
 and the density $f=F'$ is regularly varying at the infinity with index $\rho=-(1+\gamma)$, where $\gamma>0$,   moreover, we will suppose that there is no symmetry, i.e. that we are not in a~situation where $f(x)/f(-x)\to 1$, $|x|\to \infty$ and simultaneously $\an/\bn\to 1$. (If $f(x)=f(-x)$ for all sufficiently large $|x|$ and $\an=\bn$, then  $t_{j,n}=0$, $j=1,2,3$).

 Note first of all that by $[A_1]$ we have $\an=F(\qa)$, $\bn=1-F(\qb)$ and that the regularity condition implies
\begin{equation}
\label{1_gamma}
\lim_{\nty} \frac{\an}{|\qa|f(\qa)}=\lim_{\nty} \frac{\bn}{\qb f(\qb)}=\frac 1{\gamma}\, .
\end{equation}
Let $h(n)\sim g(n)$ denotes that $\lim_{\nty}h(n)/g(n)= c$, where $0<c<\infty$ is some constant.
We will now distinguish three cases: \\
(1) \ $0<\gamma<2$. In this case $\si^2(0,1) =\infty$, i.e. we are dealing with a~heavy tailed distribution $F$. Using Karamata type property (cf.~Feller~\cite{feller}, Vol.~II,
Chpt.~VIII, paragraph~9, Theorem~2), we find that  $t_{1,n}\sim\frac 1{\sqrt{n}}\frac{\an\qa^3+\bn\qb^3}{(\an\qa^2+\bn\qb^2)^{3/2}}$, the latter in absolute value is less than $\frac 1{\sqrt{n}}\left(\frac{\an|\qa|^3}{(\an\qa^2)^{3/2}}+ \frac{\bn|\qb|^3}{(\bn\qb^2)^{3/2}}\right)=\frac 1{\sqrt{k_n}}+\frac 1{\sqrt{m_n}}$. Similarly we easy check that both $t_{2,n}$ and $t_{3,n}$ are of the same order $O\left(\frac 1{\sqrt{k_n}}+\frac 1{\sqrt{m_n}}\right)$. Thus, in case (1) (in absence of symmetry) we obtain
\begin{equation}
\label{opt}
|G_n(x)-\Phi(x)|\asymp \frac{1}{\sqrt{\kn}}\ +\ \frac{1}{\sqrt{\mn}}
\end{equation}
\noindent{(2)} \ $2<\gamma<3$, the second moment of $F$ is finite, but the third moment is infinite. In this case $\si^2_{W_{(n)}}\to \si^2(0,1) <\infty$, and again using Karamata type properties of truncated moments we obtain that $t_{1,n}\sim \frac{\an\qa^3+\bn\qb^3 }{\sqrt{n}}$, and the latter quantity in absolute value is $\sim \frac 1{\sqrt{n}}\left(\an^{1-\frac 3{\gamma}}+\bn^{1-\frac 3{\gamma}}\right) L(n)$, where $L(n)$ is some positive slowly varying function. The latter quantity is the same as $n^{-3(\frac 12-\frac 1{\gamma})} \left( k_n^{1-\frac 3{\gamma}}+m_n^{1-\frac 3{\gamma}} \right)L(n) $. Note that $1-\frac 3{\gamma}<0$ since $\gamma<3$. For $t_{2,n}$ we easily find that it is of the same order as $t_{1,n}$. Note that when $\gamma$ gets close to $3$ the order of $t_{1,n}+t_{2,n}$ becomes close to $n^{-1/2}$. However, for $t_{3,n}$ we obtain a~slower rate of convergence to zero than for $t_{j,n}$, $j=1,2$. Simple computations using regularity condition show that $t_{3,n}$ in absolute value is of the order $n^{-(\frac 12-\frac 1{\gamma})}\left(k_n^{-\frac 1{\gamma}}+m_n^{-\frac 1{\gamma}}\right)L(n)$. Since $-\frac 1{\gamma}<1-\frac 3{\gamma}$ when $\gamma>2$, we see that the bias term $t_{3,n}$ is of bigger order than  $t_{1,n}+t_{2,n}$. \\
(3) \ $\gamma\ge 3$. In this case the third absolute moment of $F$ is finite and obviously $|t_{1,n}+t_{2,n}|=O(\frac 1{\sqrt{n}})$. However for the bias term  $t_{3,n}$ as before we have (in absence of symmetry) the exact order  $n^{-(\frac 12-\frac 1{\gamma})}\Big(k_n^{-\frac 1{\gamma}}+m_n^{-\frac 1{\gamma}}\Big)L(n)$, the latter quantity is close to $n^{-1/6}$ when $\gamma>3$ is close to $3$. So, in the case of a~light tailed distribution $F$ with finite third absolute moment the bias part of the Edgeworth type expansion \eqref{G_n} is of the order close to $n^{-1/6}$.

We conclude this remark by noting that in the cases (2) and (3) centering by $\boldsymbol{E}T_n$ in fact to be preferred.
However in the 'heavy tailed' case (1) this is not possible as already was shown by Peng~\cite{peng}. In a~way all this tell us that the expansion \eqref{G_n} in its present form is only really suitable for 'heavy tailed' distribution $F$, i.e. in case (1). Otherwise one should center $T_n$ by its exact expectation $\boldsymbol{E}T_n$ and consequently delete
the bias term $-\phi(x)t_{3,n}$ presented in \eqref{G_n}.

\medskip
Here is our general result on the validity of one-term Edgeworth type expansion  for the $df$'s of a normalized slightly trimmed sum.
\begin{theorem} Suppose that the conditions $[A_1]$-$[A_3]$ hold.
Then
\begin{equation}
\label{EE_n1}
\begin{split}
 \sup_{x\in \mathbb{R}}\lv F_{T_n}(x) - G_n(x) \rv
\le \frac {C_1}{n}\bigl(\delta_{1,n}&+\delta_{2,n}+\delta_{3,n}
\bigr)+ \frac {C_2}{n^{3/4}}\delta_{4,n} +\\
&+\frac {C_3}{n^{1/2}} \bigl(\delta_{5,n} +\delta_{6,n}\bigr)+
C_4\bigl(\kn^{-c}+\mn^{-c} \bigr),
\end{split}
\end{equation}
for every \ $c>0$ and some constants $C_i>0$, $i=1,\dots,4$, not depending on \ $n$, \
whereas $\delta_{1,n}=\frac{\boldsymbol{E}\lr W_1(n)\rr^4}{\si^4_{W_{(n)}}}$,
\ $\delta_{2,n}=\alp_n^2\lr \frac
1{\si_{W_{(n)}}f(\qa)}\rr^{2+\varepsilon}+\be_n^2\lr \frac
1{\si_{W_{(n)}}f(\qb)}\rr^{2+\varepsilon}$, for every $\varepsilon>0$,
while  $\delta_{3,n}$ is same as $\delta_{2,n}$, but with $\varepsilon=0$,
\ in addition  \ $\delta_{4,n}= \frac{(\ln\kn)^{5/4}
(\an)^{3/4}}{f(\qa)\si_{W_{(n)}}}+\frac{(\ln\mn)^{5/4}
(\bn)^{3/4}}{f(\qb)\si_{W_{(n)}}}$, \ \ $\delta_{5,n}=\an\frac
{\ln\kn}{\si_{W_{(n)}}}\, \Psi_{\an,\frac 1{f(x)}}(B)+ \bn\frac
{\ln\mn}{\si_{W_{(n)}}}\, \Psi_{1-\bn,\frac 1{f(x)}}(B)$, where
$B>0$ is some constant depending only on $c$.  Finally \
$\delta_{6,n}=\frac{\lv (\lambda_{1_{(n)}}+3\lambda_{2_{(n)}})\,
b_n\rv}{\si_{W_{(n)}}}$. The constants $C_i$ on the r.h.s. of \eqref{EE_n1} depend on $c$ and $\varepsilon$.
\end{theorem}

Note that  in any case we have $\frac 1n
\delta_{1,n}=O(\frac 1{\kn}+\frac 1{\mn})$ at the r.h.s. of
\eqref{EE_n1} in view of~\eqref{1.7} (cf.~proof of Theorem~1.5,
section~3).

The next corollary provides an explicit upper bound of a much simpler form.
To state it we will need the following assumption:

\medskip
\noindent$[L].$ \ {\it There exists  $0<s\le 1$ such that}
\begin{equation*}
\limsup_{\nty}\, \frac{n^s}{\kn\, \wedge\, \mn}\, \, <\, \infty \,
.
\end{equation*}
\begin{corollary} Suppose that the conditions $[A_1]$, $[A'_2]$,  and
$[L]$ hold true, in addition, assume that
 for every $0<B<\infty$ \ \ $\Psi_{\an, \frac
1{f(x)}}(B)
=O\lr\frac{1}{f(\qa)}\Bigl(\frac{\ln\kn}{\kn}\Bigr)^{1/4}\rr$  and
 $\Psi_{1-\bn, \frac 1{f(x)}}(B)
=O\lr\frac{1}{f(\qb)}\Bigl(\frac{\ln\mn}{\mn}\Bigr)^{1/4}\rr$, as $\nty$.
Then the bound on the r.h.s. in \eqref{EE_n1} is of the order $O\lr
\frac{(\ln\kn)^{5/4}}{\kn^{3/4}}+\frac{(\ln\mn)^{5/4}}{\mn^{3/4}}
\rr$.
\end{corollary}

The following theorem ensures an expansion of the Edgeworth type for  $df$ of a~normalized slightly trimmed mean.
\begin{theorem}
Suppose that $\an \vee \bn \to 0$, as $\nty$, conditions $[A_1]$
and $[L]$ hold true, and the density satisfies condition $[R]$
with $0<\gamma < 2$. Then
\begin{equation}
\label{EE_n2}
 \sup_{x\in \mathbb{R}}\lv F_{T_n}(x) - G_n(x) \rv
=O\lr
\frac{(\ln\kn)^{5/4}}{\kn^{3/4}}+\frac{(\ln\mn)^{5/4}}{\mn^{3/4}}
\rr, \quad \nty \, .
\end{equation}
\end{theorem}

The following corollary of Theorem~1.4 can be viewed as a version of the
result by Gribkova and Helmers~\cite{gh2007} under slightly weaker conditions.
\begin{corollary}
Suppose that $0<a_1<b_2<1$, where $a_1\,,b_2$ as in \eqref{seq},
condition $[A_1]$ holds true, and the density $f$ satisfies
a~H\"{o}lder condition of degree $d>0$ on the sets  $F^{-1}(U_a)$ and  $F^{-1}(U_b)$,  $U_a$ and $U_b$ as in \eqref{1.6}. Then
\begin{equation}
\label{EE_2007}
 \sup_{x\in \mathbb{R}}\lv F_{T_n}(x) - G_n(x) \rv
=o\lr n^{-1/2-p}\rr, \ \ \nty,
\end{equation}
for every $p<\min(1/4,\, d/2)$.
\end{corollary}

To proceed we state our results on the Edgeworth type expansion for a Studentized slightly trimmed sum.

Define an empirical quantile function
$\widehat{Q}_n(u)=X_{\kn:n}\vee F^{-1}_n(u)\wedge X_{n-\mn:n}$,
and the plug-in estimates of $\mu_{W_{(n)}}$ and $\si^2_{W_{(n)}}$:
\begin{equation}
\label{mu_plug} \widehat{\mu}_{W_{(n)}}=\int_0^1\widehat{Q}_n(u)\,
du=\frac{\kn}{n}X_{\kn:n}+\frac 1n\sum_{i=\kn+1}^{n-\mn}
X_{i:n}+\frac{\mn}{n} X_{n-\mn:n},
\end{equation}
\begin{equation*}
\widehat{\si}^2_{W_{(n)}}=\int_0^1(\widehat{Q}_n(u)-\widehat{\mu}_{W_{(n)}})^2\,
du=\frac{\kn}{n}X^2_{\kn:n}+\frac 1n\sum_{i=\kn+1}^{n-\mn}
X_{i:n}^2+\frac{\mn}{n} X_{n-\mn:n}^2-\widehat{\mu}_{W_{(n)}}^2.
\end{equation*}
Define the $df$ of a Studentized trimmed sum by
$F_{T_n,S}(x)=\boldsymbol{P}\Bigl(\frac{\sqrt{n}(T_n-\mu(\an,1-\bn))}{\widehat{\si}_{W_{(n)}}}
\le x \Bigr)$.
We prove that the one-term expansion for $F_{T_n,S}(x)$ is given by
\begin{equation}
\label{H_n} H_n(x)=\Phi(x)+\frac{\phi(x)}{6\sqrt{n}}
\Bigl((2x^2+1)\lambda_{1_{(n)}} +3(x^2+1)\lambda_{2_{(n)}}
-6\sqrt{n}\frac{b_n}{\si_{W_{(n)}}} \Bigr),
\end{equation}
where $b_n$ is as in \eqref{G_n}. Define a~quantity
\begin{equation}
\label{delta_s} \Delta_{n,S}=\sum_{i=1}^{5}\delta_i(n), \quad
\end{equation}
with
\begin{equation*}
\begin{split}
\delta_1(n)&=\frac{\alp_n^{3/2}}{\si_{W_{(n)}}f(\qa)}\lr\frac{\ln\kn}{\kn}\rr^{3/4}+
\frac{\be_n^{3/2}}{\si_{W_{(n)}}f(\qb)}\lr\frac{\ln\mn}{\mn}\rr^{3/4},\\
\delta_2(n) &=B\Bigl[\ln \kn\Bigl(q^2_{\an}+ \frac
{\alp_n^2}{n\si^2_{W_{(n)}}}\Psi^2_{\an,\frac 1{f(x)}}(B) \Bigr)
+\ln\mn\Bigl(q^2_{\bn}+\frac
{\be_n^2}{n\si^2_{W_{(n)}}}\Psi^2_{1-\bn,\frac 1{f(x)}}(B)
\Bigr)\Bigr],
\end{split}
\end{equation*}
where $B>0$ is some constant,
\begin{equation*}
\begin{split}
\delta_3(n)&=\frac{\alp_n^{3/2}}{\si_{W_{(n)}}}\Psi_{\an,\frac
1{f(x)}}(B)\lr\frac{\ln\kn}{\kn}\rr^{1/2}+
\frac{\be_n^{3/2}}{\si_{W_{(n)}}}\Psi_{1-\bn,\frac
1{f(x)}}(B)\lr\frac{\ln\mn}{\mn}\rr^{1/2},\\
\delta_4(n)&=\frac 1{\sqrt{n}}\lr\frac 1{\sqrt{\kn}}+\frac
1{\sqrt{\mn}}\rr\lb\frac {\alp_n^{3/2}\ln\kn}{f(\qa)
\si_{W_{(n)}}} + \frac {\be_n^{3/2}\ln\mn}{f(\qb) \si_{W_{(n)}}}
\rb,\\
\delta_5(n)&= \frac{\ln\kn \qa^2+\ln\mn
\qb^2}{n\si^2_{W_{(n)}}}=O\lr\frac {\ln\kn}{\kn}+\frac
{\ln\mn}{\mn}\rr.
\end{split}
\end{equation*}
The quantity $\Delta_{n,S}$ will determine the order of the
remainder term in the stochastic approximation for the difference
$\frac {\widehat{\si}^2_{W_{(n)}}}{\si^2_{W_{(n)}}}-1$ by a~sum of
i.i.d. of r.v.'s (cf~Lemma~2.2, Section~2):

Here is our result for a Studendized slightly trimmed sum.
\begin{theorem}
Suppose that the conditions $[A_1]$, $[A_2]$ and $[A_3]$ hold
true. Then
\begin{equation}
\label{EE_s1}
\begin{split}
 \sup_{x\in \mathbb{R}}\lv F_{T_n,S}(x) - H_n(x) \rv
\le C(\delta_{n}+\delta_{n,S}),
\end{split}
\end{equation}
where  $C>0$ is some constant not depending on $n$,\ $\delta_{n}$ is
the bound on the~r.h.s. of \eqref{EE_n1} (cf.~Theorem~1.4) and
$\delta_{n,S}=\Delta_{n,S}+\sum_{i=1}^{2}\delta_{i,S}(n)$, where
$\Delta_{n,S}$ is as in \eqref{delta_s};
$\delta_{1,S}(n)=\ln\kn(\frac 1{\sqrt{\kn}}q_{\an}+\frac 1{\kn}
)+\ln\mn(\frac 1{\sqrt{\mn}}q_{\bn}+\frac 1{\mn} )$; \
$\delta_{2,S}(n)=(\ln\kn)^2 \, q^3_{\an}+(\ln\mn)^2 \,
q^3_{\bn}+\ln\kn\ln\mn q_{\an}q_{\bn}(q_{\an}+q_{\bn})$.
\end{theorem}

The next corollary is analogous to Corollary~1.3, now for a Studentized
$T_n$.
\begin{corollary} Suppose that the conditions of Corollary~1.3 hold true. Then the bound on the r.h.s. in
\eqref{EE_s1} is of the order $O\lr
\frac{(\ln\kn)^{5/4}}{\kn^{3/4}}+\frac{(\ln\mn)^{5/4}}{\mn^{3/4}}
\rr, \ \nty\, .$
\end{corollary}

Finally, we state our Edgeworth type result for a Studentized
$T_n$ parallel to Theorem~1.5.
\begin{theorem}
Suppose that $\an \vee \bn \to 0$, as $\nty$, conditions $[A_1]$
and $[L]$ hold true, and the density satisfies condition $[R]$
with $0<\gamma < 2$. Then
\begin{equation}
\label{EE_s2}
 \sup_{x\in \mathbb{R}}\lv F_{T_n,S}(x) - H_n(x) \rv
=O\lr
\frac{(\ln\kn)^{5/4}}{\kn^{3/4}}+\frac{(\ln\mn)^{5/4}}{\mn^{3/4}}
\rr, \quad \nty\, .
\end{equation}
\end{theorem}

\medskip
 \noindent{\bf Remark~1.2} \ We conjecture that both \eqref{EE_n2} and \eqref{EE_s2}
 are also  valid without condition $[L]$. The latter condition is only used
 the formula $\delta_{2,n}$ on the r.h.s. of \eqref{EE_n1} and a~similar
  term in the studentized case.
 The \ $\varepsilon>0$ \ appearing in the  expression of estimate for
 $\delta_{2,n}$ and in its counterpart for the studentized
 case is due to the presence of a similar error term involving \
 $\varepsilon$ \
 in Bentkus et al.~\cite{bjz}. These authors, however, also conjecture
 in their Remark~1.3 that assuming the existence of  such positive \
 $\varepsilon$ \ is in fact superfluous, i.e. taking $\varepsilon =0$ will
 also work.  Without condition $[L]$ an~extra term of order
 $O\Bigl(\frac 1{\kn} \bigl(\frac{n}{\kn}\bigr)^{3\varepsilon/2}
+\frac 1{\mn} \bigl(\frac{n}{\mn}\bigr)^{3\varepsilon/2}\Bigr)$
shows up; this term can be absorbed in the r.h.s.'s of
\eqref{EE_n2} and \eqref{EE_s2} in case condition $[L]$ is
satisfied.

\medskip

  To conclude this section we want to mention a by now classical paper by van~Zwet~\cite{zwet} on Berry Esseen bounds for general symmetric statistics.  We also refer to recent work by Chen \& Shao~\cite{cs}, using a method due originally to C.Stein , and also to  Bentkus, Jing and
Zhou~\cite{bjz} who obtained optimal results on rates of convergence for $U$-statistics of general degree $k$.

For interesting recent probabilistic work on slightly trimmed sums when data are long range dependent linear processes rather than i.i.d. observations  we refer to Kulik (cf.~\cite{kulik}.

\section{A $U$-statistic approximation}
\label{U_sata}
 In this section we will approximate $T_n$ by a suitable $U$-statistic of degree~$2$. This will enable us to establish second order approximations -- \ Berry -- Esseen bounds and Edgeworth type expansions -- \ for $T_n$ and its studentized version by applying known results of this type for $U$-statistics of degree~$2$ (cf.~Friedrich~\cite{friedrich} and Bentkus et.al~\cite{bgz}). This method of proof is  well known in the literature; we refer to Bentkus et al~\cite{bjz} for recent work on this topic.  However, our remainder term -- i.e. the difference between $T_n$ and the approximating $U$-statistic -- has a~different structure compared with  the error terms appearing in previous work on 'smooth statistics' (cf.,~for instance, Putter \& van Zwet~\cite{pz}): no terms of higher order in the Hoeffding decomposition, but instead a~remainder term of Bahadur type.

Set \ $\boldsymbol{1}_{\nu}(X_i)=\boldsymbol{1}_{\{X_i\le\
\xi_{\nu}\}}$, where $\xi_{\nu}=F^{-1}(\nu)$, \ $0<\nu<1$, and
$\boldsymbol{1}_{A}$ is the indicator of the event $A$.

Define a $U$-statistic of degree $2$ with kernel, depending on
$n$, by
\begin{equation}
\label{4.1}
 L_n+U_n=\sum_{i=1}^{n}L_{n,i}+\sum_{1\, \le \, i\, <}\sum_{j\, \le \,
 n}U_{n,(i,j)}\, ,
\end{equation}
where
\begin{equation}
\label{4.2}
\begin{split}
 L_{n,i}&=\frac
1{\sqrt{n}}\bigl(W_i(n)-\mu_{W_{(n)}}\bigr) =\frac
1{\sqrt{n}}\bigl[X_i
 \boldsymbol{1}_{1-\bn}(X_i)\bigl(1-\boldsymbol{1}_{\an}(X_i)\bigr)\\
 & +\qa \boldsymbol{1}_{\an}(X_i)
 +\qb\bigl(1-\boldsymbol{1}_{1-\bn}(X_i)\bigr)- \mu_{W_{(n)}}\bigr]\, ,
\end{split}
\end{equation}
with $W_i(n)$ and $\mu_{W_{(n)}}$ as is \eqref{1.3} and
\eqref{1.5} respectively, and
\begin{equation}
\label{4.3}
\begin{split}
 U_{n,(i,j)}&=\frac
1{n\sqrt{n}}\Bigl[-\frac
1{f(\qa)}\Bigl(\boldsymbol{1}_{\an}(X_i)-\an \Bigr)
\Bigl(\boldsymbol{1}_{\an}(X_j)-\an \Bigr)\\
 & + \frac 1{f(\qb)}\Bigl(\boldsymbol{1}_{1-\bn}(X_i)-(1-\bn) \Bigr)
\Bigl(\boldsymbol{1}_{1-\bn}(X_j)-(1-\bn) \Bigr) \Bigr].
\end{split}
\end{equation}
Note that
\begin{equation}
\label{4.4}
\boldsymbol{E} L_{n,i}=0,\ \ i=1,\dots,n\, ,
\end{equation}
and
\begin{equation}
\label{4.5} \boldsymbol{E} U_{n,(i,j)}=0,\ \ \boldsymbol{E}\bigl(L_{n,i}
U_{n,(i,j)}\bigr)=0, \ \ i,j=1,\dots,n \ \,(i\neq j\, )\, .
\end{equation}

Using \eqref{4.1}--\eqref{4.5}, we easily check that $\boldsymbol{E}
\bigl(L_n+U_n\bigr)^2=\si^2_{W_{(n)}}+ \boldsymbol{E} (U^2_n)$, where
$\si^2_{W_{(n)}}$ is given as in \eqref{1.5} and
$\boldsymbol{E}(U^2_n)=\frac{n-1}{2n^2}\, \boldsymbol{E}\bigl(
n^{3/2}U_{n,(1,2)} \bigr)^2\le \frac 1n \lr
\frac{\alp_n^2}{f^2(\qa)} +\frac{\be_n^2}{f^2{(\qb)}}\rr$. So we
obtain:
\begin{equation}
\label{4.6} \boldsymbol{E} \lr \frac{L_n + U_n}{\si_{W_{(n)}}
}\rr^2=1+\varepsilon_n\, ,
\end{equation}
where $0<\varepsilon_n\le q^2_{\an}+q^2_{\bn}$ (cf.~\eqref{qu}), and
$\varepsilon_n\to 0$, as $\nty$, provided condition $[A_2]$ is
satisfied.

For the third moment we have
\begin{equation}
\label{4.7} \boldsymbol{E} \lr \frac{L_n + U_n}{\si_{W_{(n)}} }\rr^3
=\frac 1{\sqrt{n}}\lambda_{1,n}+3\si^{-3}_{W_{(n)}}\,
\boldsymbol{E}(L^2_n U_n) +3\si^{-3}_{W_{(n)}}\, \boldsymbol{E}
(L_nU^2_n)+\si^{-3}_{W_{(n)}}\, \boldsymbol{E} (U^3_n),
\end{equation}
where $\lambda_{1,n}$ is as in \eqref{la_12}. For the second term
on the r.h.s. of \eqref{4.7} we obtain $3\si^{-3}_{W_{(n)}}\,
\boldsymbol{E}(L^2_n U_n)=3n(n-~1)
\si^{-3}_{W_{(n)}}\boldsymbol{E}(L_{n,1}L_{n,2}U_{n,(1,2)})$, which is
equal to
\begin{equation}
\label{4.8} \begin{split} & 3\, {\si^{-3}_{W_{(n)}}}\, \frac
{n-1}{n\sqrt{n}}\lb -\frac{(\an
(\qa-\mu_{W_{(n)}}))^2}{f(\qa)}+\frac{(\bn
(\qb-\mu_{W_{(n)}}))^2}{f(\qb)} \rb\\
=\, & 3\, \frac 1{\sqrt{n}}\lambda_{2,n}\lr1-\frac 1n \rr=3\frac
1{\sqrt{n}}\lambda_{2,n} +o\lr\frac 1n\rr,
\end{split}
\end{equation}
where $\lambda_{2,n}$ is as in \eqref{la_12}. The last equality
on the r.h.s. of \eqref{4.8} is valid because by \eqref{1.7} there
exists a constant $C>0$ such that $\frac 1{n\sqrt{n}}\lv
\lambda_{2,n}\rv \le \frac
Cn\lr\frac{\an}{\sqrt{n}\si_{W_{(n)}}f(\qa)
}+\frac{\bn}{\sqrt{n}\si_{W_{(n)}}f(\qb) }\rr=O\lr\frac
1n(q_{\an}+q_{\bn})\rr=o\lr\frac 1n\rr$. Using relations
\eqref{4.2}--\eqref{4.5}  we find that $ \frac
3{\si^3_{W_{(n)}}}\, \boldsymbol{E}\bigl(L_n U^2_n\bigr) =\frac
{3(n-1)}{\si^3_{W_{(n)}}n^2\sqrt{n}}\lb \frac{\alp_n^2\,
(\qa-\mu_{W_{(n)}})\, (1-\an)\, (1-2\an)}{f^2(\qa)}
-\frac{2\alp_n^2\be_n^2\,
\bigl[(\qa-\mu_{W_{(n)}})+(\qb-\mu_{W_{(n)}})\bigr]
}{f(\qa)f(\qb)}\rp$ $\lp+\frac{\be_n^2\, (\qb-\mu_{W_{(n)}})\,
(1-\bn)\, (1-2\bn)}{f^2(\qb)}\rb$, and by relation \eqref{1.7} the
latter quantity is of the order $O\lr q^2_{\an}\frac
1{\sqrt{\kn}}+ q^2_{\bn}\frac
1{\sqrt{\mn}}+q_{\an}q_{\bn}\an\bn\bigl( \frac 1{\sqrt{\kn}}+\frac
1{\sqrt{\mn}}\bigr)\rr$. Finally, for the fourth term at the
r.h.s. of \eqref{4.7} we have ${\si^{-3}_{W_{(n)}}}\,
\boldsymbol{E}\bigl(U^3_n\bigr)=\si^{-3}_{W_{(n)}}\frac {(n-1)}{2\,
n^3\sqrt{n}}\boldsymbol{E}\bigl(n\sqrt{n}\, U^3_{n,(1,2)}\bigr)$, and
after simple computations we obtain that the latter quantity in
absolute value is less than $\si^{-3}_{W_{(n)}}\frac {1}{2\,
n^2\sqrt{n}}\lb\frac {\alp_n^2}{f^3(\qa)}+ \frac
{\be_n^2}{f^3(\qb)}+ 3 \frac
{\alp_n^2\be_n^2}{f^2(\qa)f(\qb)}+3 \frac
{\alp_n^2\be_n^2}{f(\qa)f^2(\qb)}\rb=r_1+r_2$, where $r_1=\frac
12\lr\frac 1{\kn}q^3_{\an}+\frac 1{\mn}q^3_{\bn}\rr$, and
$r_2=\frac 3{2\, n} q_{\an}\, q_{\bn}\lr\bn\, q_{\an}+ \an\,
q_{\bn}\rr$. Hence, under conditions $[A_1]$--$[A_2]$ we have
$r_1=o\lr q^2_{\an}\frac 1{\sqrt{\kn}}+ q^2_{\bn}\frac
1{\sqrt{\mn}}\rr$ and $r_2=o\lr \frac 1n (q_{\an}+ q_{\bn})\rr =o\lr \frac 1n\rr$,
and we can conclude that the  fourth term on
the r.h.s. of \eqref{4.7} is of negligible order for our purposes.

Our computations directly imply that
\begin{equation}
\label{4.9} \boldsymbol{E} \lr \frac{L_n + U_n}{\si_{W_{(n)}}
}\rr^3=\frac 1{\sqrt{n}}\bigl(
\lambda_{1,n}+3\lambda_{2,n}\bigr)+R_n\, ,
\end{equation}
where $R_n=O\lr q^2_{\an}\frac 1{\sqrt{\kn}}+ q^2_{\bn}\frac
1{\sqrt{\mn}}+ \an\, \bn \, q_{\an}\, q_{\bn}\bigl(\frac
1{\sqrt{\kn}}+\frac 1{\sqrt{\mn}} \bigr)\rr$.

The next lemma provides an estimate of the precision of the
approximation of $T_n$ by the sum of a~$U-$statistic with varying kernel of
the form \eqref{4.1} with mean zero and a~bias term $b_n$.
\begin{lemma}
Suppose that the conditions $[A_1]$ and $[A_2]$ hold true. Then
\begin{equation}
\label{L4.1_1} \boldsymbol{P}\Bigl(\big|n^{1/2}\bigl(T_n-\mu(\an,1-\bn)\bigr)
-(L_n+U_n+b_n)\big|>\Delta_n \Bigr)=O\bigl((\kn
\wedge\mn)^{-c}\bigr),
\end{equation}
for every $c>0$, where $b_n$ is as in \eqref{G_n}, \
$\Delta_n=A(\Delta_{\alp,n}+\Delta_{\be,n})$,
\begin{equation}
\label{L4.1_2}
\begin{split}
\Delta_{\alp,n}=&\, \an\, \frac{\ln\kn}{\sqrt{n}}\Bigl[\frac
1{f(\qa)}\lr\frac {\ln\kn}{\kn}\rr^{1/4} +\Psi_{\an,\frac 1{f(x)}}(B)\Bigr],\\
\Delta_{\be,n}=&\, \bn\, \frac{\ln\mn}{\sqrt{n}}\Bigl[\frac
1{f(\qb)}\lr\frac {\ln\mn}{\mn}\rr^{1/4} +\Psi_{1-\bn,\frac
1{f(x)}}(B)\Bigr],
\end{split}
\end{equation}
and where the constants $A,B>0$ depend only on $c$.
\end{lemma}
\noindent{\bf Proof.} Define a binomial r.v. $N_\nu= \sharp \{i :
X_i \le \xi_{\nu}\}$, $0<\nu<1$, and note that
\begin{equation}
\label{L4.1_3} W_n=\frac 1n\sum_{i=1}^{n}W_i(n)=\frac
{N_{\an}}{n}\qa+\frac
1n\sum_{i=N_{\an}+1}^{N_{1-\bn}}X_{i:n}+\frac{n-N_{1-\bn}}{n}\qb.
\end{equation}
Then
\begin{equation}
\label{L4.1_4}
\begin{split}
&T_n-\mu(\an,1-\bn) -[W_n-\boldsymbol{E}W_n]\\
=&\frac 1n\Bigl[ sgn(N_{\an}-\, \kn)\sum_{i=(\kn\wedge
N_{\an})+1}^{N_{\an}\vee
\kn}(X_{i:n}-\qa)\\
&-sgn(N_{1-\bn}-(n-\mn))\sum_{i=((n-\mn)\wedge
N_{1-\bn})+1}^{N_{1-\bn}\vee \, (n-\mn)}(X_{i:n}-\qb)\Bigr],
\end{split}
\end{equation}
where $sgn(s)=s/|s|$, $sgn(0)=0$, and by lemma~\ref{lem_A.2} (cf. Section~\ref{Bahadur} ) the latter is
equal to
\begin{equation}
\label{L4.1_5} -\frac{(N_{\an}-\an\, n)^2}{2\, n^2}\frac
1{f(\qa)}+\frac{(N_{1-\bn}-(1-\bn)\, n)^2}{2\, n^2}\frac
1{f(\qb)}+R_n,
\end{equation}
where $\boldsymbol{P}\Bigl(|R_n|>\frac
A{\sqrt{n}}(\Delta_{\alp,n}+\Delta_{\be,n})\Bigr)=O\lr\lr\kn\wedge
\mn\rr^{-c}\rr$, and $\Delta_{\alp,n}$, $\Delta_{\be,n}$ are given
as in \eqref{L4.1_2}. Relations \eqref{4.1}--\eqref{4.3}, and
\eqref{L4.1_4}--\eqref{L4.1_5} yield
\begin{equation}
\label{L4.1_6} \begin{split} & n^{1/2}\bigl(T_n
-\mu(\an,1-\bn)\bigr)=L_n+U_n-\frac 1{2n\sqrt{n}}\lb\frac
1{f(\qa)}\sum_{i=1}^{n}\bigl(\boldsymbol{1}_{\an}(X_i)-\an
\bigr)^2 \rp\\
\ &+\lp \frac
1{f(\qb)}\sum_{i=1}^{n}\bigl(\boldsymbol{1}_{1-\bn}(X_i)-(1-\bn)
\bigr)^2 \rb+n^{1/2}R_n\\
= \ &L_n+U_n+b_n+\frac 1{2\, \sqrt{n}}\,
\overline{r}_{n}+n^{1/2}R_n ,
\end{split}
\end{equation}
where $b_n$ is as in \eqref{G_n} \ and \ $\overline{r}_{n}=-\frac
1{f(\qa)}\overline{r}_{n,1}+\frac 1{f(\qb)}\overline{r}_{n,2}$,
\begin{equation*}
 \begin{split}
\overline{r}_{n,1}=&\frac
1n\sum_{i=1}^{n}\Bigl[\bigl(\boldsymbol{1}_{\an}(X_i)-\an \bigr)^2
- \an\, (1-\an)\Bigr],\\
\overline{r}_{n,2}=&\frac
1n\sum_{i=1}^{n}\Bigl[\bigl(\boldsymbol{1}_{1-\bn}(X_i)-(1-\bn)
\bigr)^2 - \bn\, (1-\bn)\Bigr]\, .
\end{split}
\end{equation*}
We consider only  $\overline{r}_{n,1}$, the treatment for
$\overline{r}_{n,2}$ is similar. Note that $\overline{r}_{n,1}$ is
an average of i.i.d. centered r.v.'s, $\overline{r}_{n,1}=\frac 1n
S_{n,1}$, where $S_{n,1}=\sum_{k=1}^nY_k$, \ $\boldsymbol{E}Y_k=0$,
and $B_n=D(S_{n,1})=n\si^2_1$ with
$\si^2_1=\boldsymbol{E}Y^2_1=\an\bigl(1-\an\bigr)\bigl[1-2\an\bigr]^2$.
Moreover, for each integer $m\ge 2$ we have
\begin{equation*}
\boldsymbol{E}Y_1^m=\si^2_1\bigl[1-2\, \an\bigr]^{m-2}\bigl[
\bigl(1-\an\bigr)^{m-1}+(-1)^m (\an)^{m-1}\bigr],
\end{equation*}
and hence, $\lv\, \boldsymbol{E}Y_1^m\rv\le\si^2_1$. Then by
applying an~exponential bound (cf.~Petrov~\cite{petrov},
chapter~3,Theorem~17, with $H=1$) we obtain
\begin{equation}
\label{L4.1_7} \boldsymbol{P}\Bigl(\lv S_{n,1}\rv \ge x\Bigr)\le
\exp\lr-\frac{x^2}{4B_n} \rr
\end{equation}
for every $0\le x\le B_n$. Take $x=A\bigl(n\ln\kn\,
\an(1-\an)\bigr)^{1/2}\big|1-2\, \an\big|$. If $\alp=1/2$ is not
a~partial limit point of the sequence $\an$, we can easily see
that $0\le x\le B_n$ for all sufficiently large $n$. Otherwise
note that we can consider only cases when $\delta_n=\big|\an-\frac
12\big|>A_1\lr\frac{\ln\, \kn}{\kn}\rr^{1/2}$ with some $A_1>0$
which we will choose later. Indeed, if it is not so, we can write:
\ $\overline{r}_{n,1}=\frac 1n\sum_{i=1}^{n}\Bigl[(\frac 12
+\delta_{n,i})^2 -(\frac 12-\delta_n)(\frac 12+\delta_n)\Bigr]$,
where $\delta_{n,i}=(-1)^{\boldsymbol{1}_{\an}(X_i)}\lr \an -\frac
12 \rr$, \ $|\delta_{n,i}|=\delta_n$, and $|\overline{r}_{n,1}|\le
\delta_n(1+2\,\delta_n)=O\lr\frac{\ln\kn}{\kn}\rr^{1/2}$. Since
$\xi_{1/2}\in U_a\,\cup\, U_b$ under condition $[A_1]$, we have
$f(x)\ge f_0>0$ in some neighbourhood of $\xi_{1/2}$\, , and  we
obtain $\overline{r}_{n,1}\frac
1{f(\qa)}=O\lr\bigl(\frac{\ln\kn}{\kn}\bigr)^{1/2}\rr$, and hence
$\frac 1{2\sqrt{n}}|\overline{r}_{n,1}|=o(\Delta_n)$
(cf.~\eqref{L4.1_6}).

For the case $\delta_n>A_1\lr\frac{\ln\, \kn}{\kn}\rr^{1/2}$ we
easily check that $x<B_n$ for all sufficiently large $n$ when we
choose $A_1$ such that $A_1>A(1-a_2)^{-1/2}$, and by
\eqref{L4.1_7} we obtain $ \boldsymbol{P}\Bigl(\lv \overline{r}_{n,1} \rv \ge
A\frac 1{\sqrt{n}}|1-2\an|
 \bigl(\ln\kn\an(1-\an)\bigr)^{1/2} \Bigr)\le \exp\lr-\frac
 14A^2\ln\kn\rr$, the latter is of the order $O\lr\kn^{-c}\rr$
 when $A^2>4c$. So $\frac 1{f(\qa)}\lv\overline{r}_{n,1}\rv\le
 A\frac{(\ln\kn)^{1/2}\an}{\sqrt{n}f(\qa)(\an)^{1/2}}=
 A\, \an\lr\frac{\ln\kn}{\kn}\rr^{1/2}$, and hence $\frac 1{2\sqrt{n}}
 \frac 1{f(\qa)}\overline{r}_{n,1}=o(\Delta_n)$ (cf. \eqref{L4.1_2}
 and \eqref{L4.1_6}). The lemma is proved. $\quad \square$

\bigskip
We complete this section by a lemma which can be viewed as
an~extension of Lemma~5.1 from Gribkova and Helmers~\cite{gh2006} to
slightly trimmed means, i.e.~to the case corresponding to the
first two lines of~\eqref{1.6}.

\begin{lemma} Suppose that the conditions $[A_1]$ and $[A_2]$ hold true .
Then for every $c>0$
\begin{equation}
\label{lem_23}
 \boldsymbol{P}\Bigl(\Big|\frac
{\widehat{\si}^2_{W_{(n)}}}{\si^2_{W_{(n)}}}-1-\frac
{V_n}{\si^2_{W_{(n)}}}\Big|>A\Delta_{n,S} \Bigr)=O\Bigl(\kn^{-c}+
\mn^{-c}\Bigr),
\end{equation}
where $\Delta_{n,S}$ is as in \eqref{delta_s},
\begin{equation}
\label{V_n} V_n=V_{n,1}+V_{n,2},
\end{equation}
with
\begin{equation*}
\begin{split}
V_{n,1}=2\frac {\an}{f(\qa)}\frac{N_{\an}-\an
n}{n}\bigl(\mu_{W_{(n)}}-\qa\bigr)  \qquad \qquad \qquad \qquad \qquad \qquad \qquad \qquad\\
\qquad \qquad \qquad \qquad \qquad \qquad  +\, 2\, \frac
{\bn}{f(\qb)}\frac{N_{1-\bn}-(1-\bn)
n}{n}\bigl(\mu_{W_{(n)}}-\qb\bigr)
\end{split}
\end{equation*}
and
\begin{equation*}
V_{n,2}=\frac 1n\sum_{i=1}^{n}\lb\bigl(W_i(n)-\mu_{W_{(n)}}
\bigr)^2-\si^2_{W_{(n)}}\rb,
\end{equation*}
$A,\,B>0$ ($B$ a~constant appearing in $\delta_i(n)$, $i=2,3$) are
some constants not depending on $n$. Moreover,
\begin{equation}
\label{EV_n} \boldsymbol{E}V_n=0,\qquad \ \ \boldsymbol{E}\Bigl(\frac
{V_n}{\si^2_{W_{(n)}}}\Bigr)^2=O\Bigl(
\frac{\boldsymbol{E}(W_1(n))^4}{n\,\si^4_{W_{(n)}}}+q^2_{\an}+q^2_{\bn}\Bigr).
\end{equation}
\end{lemma}
\noindent{\bf Proof.} First we note that relations
\eqref{EV_n} follow directly by definition \eqref{V_n} of
$V_{n,i}$, $i=1,2$, and \eqref{1.7}.

To prove \eqref{lem_23} fix an arbitrary $c>0$ and define the
auxiliary quantity $S^2_{W_{(n)}}= \frac
1n\sum_{i=1}^{n}W^2_i(n)-\overline{W}^{\, 2}(n)$, where
$\overline{W}(n)=\frac 1n \sum_{}^{}W_i(n)$. First we prove that
\begin{equation}
\label{lem23_21}
\widehat{\si}^2_{W_{(n)}}=S^2_{W_{(n)}}+V_{n,1}+R_{n,1},
\end{equation}
where
\begin{equation}
\label{lem23_22}
\boldsymbol{P}\Bigl(\frac{|R_{n,1}|}{\si^2_{W_{(n)}}}>A_1\Delta_{n,S}\Bigr)=O(\kn^{-c}+\mn^{-c}).
\end{equation}
Here and elsewhere $A,A_i>0$, $i=1,2,\dots$, denote the constants,
independent of $n$. We have
\begin{equation}
\label{lem23_23}
\begin{split}
\widehat{\si}^2_{W_{(n)}}-S^2_{W_{(n)}}= \lb \frac {\kn}n
X^2_{\kn:n} + \frac 1n \sum_{i=\kn+1}^{n-\mn} X^2_{i:n} + \frac
{\mn}n X^2_{n-\mn:n}  - \frac {N_{\an}}n \qa^2  -\frac 1n
\sum_{i=N_{\an}+1}^ {N_{1-\bn}} X^2_{i:n} \rp \qquad \qquad \qquad \qquad \\
\lp\phantom{\sum_{i_1}^{j^1}}-\frac {n-N_{1-\bn}}n \qb^2 \rb +\lb
\lr \frac {N_{\an}}n \qa +\frac 1n \sum_{i=N_{\an}+1}^ {N_{1-\bn}}
X_{i:n} + \frac
{n-N_{1-\bn}}n \qb \rr^2\rp\qquad \qquad \qquad \qquad \\
\qquad \qquad \qquad \qquad \qquad\lp\phantom{\sum_{i_1}^{j^1}}
-\lr\frac {\kn}n X_{\kn:n} + \frac 1n \sum_{i=\kn+1}^{n-\mn}
X_{i:n} + \frac {\mn}n X_{n-\mn:n}\rr^2 \rb .
\qquad \qquad \qquad \qquad\ 
\end{split}\end{equation}
Rewrite the term within the first square brackets on the r.h.s. of
\eqref{lem23_23} as
\begin{equation}
\label{lem23_24}
\begin{split}
 &\frac {\kn}n \bigl(X^2_{\kn:n} - \qa^2\bigr) + \frac 1n sgn(N_{\an}-\kn)
 \sum_{i=(\kn\wedge N_{\an})+1}^{\kn \vee N_{\an}}
\bigl(X^2_{i:n}-\qa^2\bigr)\, +\\
+&\frac {\mn}n \bigl(X^2_{n-\mn:n} - \qb^2\bigr) \\
-& \frac 1n sgn(N_{1-\bn}-(n-\mn))
 \sum_{i=(N_{1-\bn}\wedge (n-\mn))+1}^{(n-\mn) \vee N_{1-\bn}}
\bigl(X^2_{i:n}-\qb^2\bigr),
\end{split}\end{equation}
then by Lemmas~\ref{lem_A.1} and ~\ref{lem_A.2}, where $G(x)=x^2$ (cf.~Section~\ref{Bahadur} ), the
latter quantity is equal to
\begin{equation}
\label{lem23_25}
\begin{split}
 &-2 \an\frac {N_{\an}-\an n}{n} \frac {\qa}{f(\qa)}-\frac {(N_{\an}-\an
 n)^2}{n^2} \frac {\qa}{f(\qa)}\\
&-2 \bn \frac {N_{1-\bn}-(1-\bn) n}{n}\frac {\qb}{f(\qb)}\\
&+\frac {(N_{1-\bn}-(1-\bn)
 n)^2}{n^2} \frac {\qb}{f(\qb)} +R^{(1)}_{n,1},
\end{split}\end{equation}
where $R^{(1)}_{n,1}$ is a remainder term appearing as result of
application of Lemma~~\ref{lem_A.1} two times: in the first and third terms
of \eqref{lem23_24}. Using \eqref{(A.1)}, \eqref{1.7} and inequalities
\eqref{eq_psi} we obtain that
$|R^{(1)}_{n,1}|/\si^2_{W_{(n)}}=O\bigl(\delta_1(n)+\delta_2(n)+\delta_3(n)
\bigr)$ with probability $1-O(\kn^{-c}+\mn^{-c})$, where
$\delta_2(n),\ \delta_3(n)$ involve $B>0$, which depends on $c$
and does not depend on $n$. Note that the remainder term appearing
as result of application of Lemma~\ref{lem_A.2} in \eqref{lem23_24} is of
the negligible order and contribute to $R^{(1)}_{n,1}$. Moreover,
the Bernstein's inequality and \eqref{1.7} together imply that
$\frac {(N_{\an}-\an
 n)^2}{n^2\si^2_{W_{(n)}}}|\qa| \frac 1{f(\qa)}=O\lr
 \frac {\alp_n^{1/2}\ln\kn}{n\si_{W_{(n)}}}\frac 1{f(\qa)}\rr$
 with probability $1-O(\kn^{-c})$, and
 it is $o\bigl(\delta_1(n) \bigr)$. The same is valid for the
 second quadratic term in \eqref{lem23_25}. Thus, we obtain that
 \eqref{lem23_25}
 is equal
\begin{equation}
\label{lem23_26}
 -2 \an\frac {N_{\an}-\an n}{n} \frac {\qa}{f(\qa)}
-2 \bn \frac {N_{1-\bn}-(1-\bn) n}{n}\frac {\qb}{f(\qb)}
+R^{(1)}_{n,1}.
\end{equation}

Now consider the term within the second square brackets on the
r.h.s. of  \eqref{lem23_23}. Arguing as before, we can rewrite it
as
\begin{equation}
\label{lem23_27}
\begin{split}
 &\lr \frac{\an}{f(\qa)}\frac {N_{\an}-\an n}{n}
 + \frac{\bn}{f(\qb)} \frac {N_{1-\bn}-(1-\bn) n}{n}
+R^{(2)}_{n,1}\rr\times\\
\times&\lr \frac 2n\sum_{i=1}^{n}W_i(n)-\frac{\an}{f(\qa)}\frac
{N_{\an}-\an n}{n} - \frac{\bn}{f(\qb)} \frac {N_{1-\bn}-(1-\bn)
n}{n} -R^{(2)}_{n,1}\rr,
\end{split}
\end{equation}
where by Lemma~~\ref{lem_A.1}  $|R^{(2)}_{n,1}|\le A_1\lf\alp_n^2\lb\frac
1{f(\qa)}\lr\frac{\ln\kn}{\kn}\rr^{3/4}+\Psi_{\an,\frac
1{f(x)}}(B)\lr\frac{\ln\kn}{\kn}\rr^{1/2} \rb \rp$ $\lp+
\be_n^2\lb\frac
1{f(\qb)}\lr\frac{\ln\mn}{\mn}\rr^{3/4}+\Psi_{\bn,\frac
1{f(x)}}(B)\lr\frac{\ln\mn}{\mn}\rr^{1/2} \rb\rf$ with
probability $1-O\bigl(\kn^{-c}+\mn^{-c}\bigr)$. The quadratic and
remainder terms caused by application of Lemma~\ref{lem_A.2} are of the
negligible order and contribute to $R^{(2)}_{n,1}$ again. Simple
computations using \eqref{lem23_23}, \eqref{lem23_26} --
\eqref{lem23_27}, Bernstein's inequality and \eqref{1.7} lead to
the following relation
\begin{equation}
\label{lem23_28}
\begin{split}
 \frac {\widehat{\si}^{\,2}_{W_{(n)}}}{\si^2_{W_{(n)}}}-\frac{S^2_{W_{(n)}}}
 {\si^2_{W_{(n)}}}
 =\frac{V_{n,1}}{\si^2_{W_{(n)}}}
 +O\Bigl(\delta_1(n)+\delta_2(n)+\delta_3(n)\Bigr) +\frac{R^{(3)}_{n,1}}{\si^2_{W_{(n)}}},
\end{split}
\end{equation}
where $R^{(3)}_{n,1}=\frac
2n\sum_{i=1}^{n}\bigl(W_i(n)-\mu_{W_{(n)}}\bigr)
\lb\frac{\an}{f(\qa)}\frac {N_{\an}-\an n}{n}
 + \frac{\bn}{f(\qb)} \frac {N_{1-\bn}-(1-\bn) n}{n}\rb$. Since by
 Hoeffding's inequality for sum of i.i.d. centered bounded r.v.'s
 (cf.~Hoeffding~\cite{ho})
 $\lv \sum_{i=1}^{n}\bigl(W_i(n)-\mu_{W_{(n)}}\bigr)\rv\le A_2\ln^{1/2}(\kn\wedge\mn)
 \bigl(\qb-\qa \bigr)n^{1/2}$ with probability
 $1-O\bigl(\kn^{-c}+\mn^{-c}\bigr)$,where $A_2>0$ is some constant,
 depending on $c$ and not depending on $n$, using the latter bound and
 Bernstein's inequality and \eqref{1.7} after the simple computations we obtain
\begin{equation}
\label{lem23_29}
 \lv\frac{R^{(3)}_{n,1}}{\si^2_{W_{(n)}}}\rv \le A_2\frac {\lv\qa\rv+\lv\qb\rv}{n\, \si^2_{W_{(n)}}} \lb \frac{\alp_n^{3/2}\, \ln\,\kn}{f(\qa)}
 + \frac{\be_n^{3/2}\, \ln\,\mn}{f(\qb)}\rb
 =O\Bigl(\delta_4(n)\Bigr)
\end{equation}
with probability $1-O\bigl(\kn^{-c}+\mn^{-c}\bigr)$, and
\eqref{lem23_21}--\eqref{lem23_22} follow.

Finally we prove that
\begin{equation}
\label{lem23_30} S^2_{W_{(n)}}=\si^2_{W_{(n)}}+V_{n,2} +R_{n,2},
\end{equation}
where $R_{n,2}$ satisfies \eqref{lem23_22}. We have
\begin{equation}
\label{lem23_31} S^2_{W_{(n)}}-\si^2_{W_{(n)}}-V_{n,2}=
S^2_{W_{(n)}}-\frac 1n
\sum_{i=1}^{n}\bigl(W_i(n)-\mu_{W_{(n)}}\bigr)^2=-\bigl(\overline{W}(n)-\mu_{W_{(n)}}\bigr)^2,
\end{equation}
and applying Hoeffding's inequality once more, we obtain that the
quantity at the r.h.p. of \eqref{lem23_31} divided by
$\si^2_{W_{(n)}}$ in absolute value is of the order
$\frac{\ln\kn\wedge\ln\mn}{n\,
\si^2_{W_{(n)}}}\Bigl(\qb-\qa\Bigr)^2= O\Bigl(\delta_5(n)\Bigr)$
with probability $1-O\bigl(\kn^{-c}+\mn^{-c}\bigr)$, and
\eqref{lem23_30} follows. The lemma is proved.$\quad \square$
\section{Proofs}
\label{Proofs}
In this section we prove Theorems 1.1.-1.7 and their corollaries
stated in Section 1.

\medskip
\noindent {\bf Proof of Theorem 1.1.} By Lemma~2.1 we can write
$n^{1/2}(T_n-\mu(\an,1-\bn))=L_n+U_n+b_n+R_n$, where $L_n+U_n$ is
$U$-statistic of degree 2, (cf.~\eqref{4.1}), $b_n$ is as in
\eqref{G_n}, and $R_n$ is a remainder term (cf.~\eqref{L4.1_1}).
Define the $df$ of a~normalized $U$-statistic:
$F_{U,n}(x)=\boldsymbol{P}\lr\frac{L_n+U_n}{\si_{W_{(n)}}}\le x\rr$.

Since $\frac{|b_n|}{\si_{W_{(n)}}}\le\frac 12( q_{\an}+q_{\bn})=\frac
1{2\sqrt{n}}\delta_{2,n}$ \ (cf.~\eqref{BE_1}), the following
inequalities are valid:
\begin{equation}
\label{petr_l}
 F_{U,n}(x-\delta_n)-\boldsymbol{P}\lr |R_n|> \Delta_n \rr\le
F_{T_n}(x)\le F_{U,n}(x+\delta_n)+\boldsymbol{P}\lr |R_n|> \Delta_n \rr,
\end{equation}
where $\delta_n=\frac{\Delta_n}{\si_{W_{(n)}}}+\frac
1{2\sqrt{n}}\delta_{2,n}$, \ $\Delta_n$ is as in \eqref{L4.1_1},
and by Lemma~2.1 $\boldsymbol{P}\lr |R_n|> \Delta_n
\rr=O\lr\kn^{-c}+\mn^{-c}\rr$ for every $c>0$.

For $F_{U,n}(x\pm\delta_n)$ we can write
\begin{equation}
\label{fi1} \sup_{x\in\mathbb{R}}\lv
F_{U,n}(x\pm\delta_n)-\Phi(x)\rv \le \Delta_{n,1}+\Delta_{n,2},
\end{equation}
where
\begin{equation}
\label{fi2} \Delta_{n,1}=\sup_{x\in\mathbb{R}}\lv
F_{U,n}(x)-\Phi(x)\rv ,\quad \Delta_{n,2}=\sup_{x\in\mathbb{R}}\lv
\Phi(x\pm\delta_n)-\Phi(x)\rv.
\end{equation}
To estimate $\Delta_{n,1}$ we apply the~Berry -- Esseen bound for
$U$-statistics (cf.~Friedrich~\cite{friedrich}):
\begin{equation}
\label{fi3} \Delta_{n,1}\le \frac
C{\sqrt{n}}\lb\frac{\boldsymbol{E}\lv W_1(n) \rv^3}{\si^3_{W_{(n)}}}+
\frac {\boldsymbol{E}\lv
n\sqrt{n}U_{n,(1,2)}\rv^{5/3}}{\si^{5/3}_{W_{(n)}}}\rb,
\end{equation}
where $C>0$ is an absolute constant. Using formula \eqref{4.3}, we
easily check that
\begin{equation}
\label{fi4} \begin{split}\boldsymbol{E} \lv
n\sqrt{n}U_{n,(1,2)}\rv^{5/3}\le & \, 2^{2/3}\lb\lr\frac
1{f(\qa)}\rr^{5/3}\Bigl(\boldsymbol{E}\lv\boldsymbol{1}_{\an}(X_1)-\an\rv^{5/3}
\Bigl)^2\rp \\
&\lp + \lr\frac
1{f(\qb)}\rr^{5/3}\Bigl(\boldsymbol{E}\lv\boldsymbol{1}_{1-\bn}(X_1)-(1-\bn)\rv^{5/3}
\Bigl)^2\rb\\
\le &\, 2^{4/3} \lb \frac{\bigl(\an
(1-\an)\bigr)^2}{f^{5/3}(\qa)}+ \frac{\bigl(\bn
(1-\bn)\bigr)^2}{f^{5/3}(\qb)}\rb.
\end{split}
\end{equation}
Relations \eqref{fi3}-\eqref{fi4} together imply that
\begin{equation}
\label{fi5} \Delta_{n,1}\le \frac
{C_1}{\sqrt{n}}\bigl(\delta_{1,n} +\delta_{3,n}\bigr),
\end{equation}
where $C_1>0$ is some absolute constant.

Finally, consider $\Delta_{n,2}$. Note that
$\frac{\Delta_n}{\si_{W_{(n)}}}=\frac 1{\sqrt{n}}\bigl(o\lr
\delta_{2,n} \rr+\delta_{4,n}\bigr)$ (cf.~\eqref{BE_1}
and~\eqref{L4.1_1}), therefore $\delta_n=\frac
1{\sqrt{n}}\bigl(\delta_{2,n}(\frac 12+o(1)) +\delta_{4,n}\bigr)$,
and we obtain
\begin{equation}
\label{fi6} \Delta_{n,2}\le \frac
{C_2}{\sqrt{n}}\bigl(\delta_{2,n} +\delta_{4,n}\bigr),
\end{equation}
where $C_2>0$ is some  constant, depending only on $c$
(cf.~Lemma~2.1). Relations \eqref{petr_l}--\eqref{fi1}, and
\eqref{fi5}--\eqref{fi6} imply \eqref{BE_1}. The theorem is
proved.$\quad \square$

\noindent {\bf Proof of Theorem 1.2.} We obtain this theorem as a
consequence of Theorem~1.1. First choose $c=1/2$ in \eqref{BE_1}. To
prove \eqref{BE_2} we must verify that under conditions $[A_1]$,
$[A'_2]$ and $[A'_3]$
\begin{equation*}
\label{pt2_1} \frac 1{\sqrt{n}}\, \delta_{i,n}=O\Bigl(\frac
1{\sqrt{\kn}}+\frac 1{\sqrt{\mn}}\Bigr) ,\quad i=1,\dots,4.
\end{equation*}
For $i=1$ we have
\begin{equation}
\label{pt2_2}  \frac 1{\sqrt{n}}\, \delta_{1,n}=\frac{\an
|\qa|^3+\int_{\an}^{1-\bn}|F^{-1}(u)|^3\,du +\bn
|\qb|^3}{\sqrt{n}\si^3_{W_{(n)}}}.
\end{equation}
We consider the three terms in the nominator on the r.h.s. of
\eqref{pt2_2}. For the first one we have $\frac{\an
|\qa|^3}{\sqrt{n}\, \si^3_{W_{(n)}}}=\frac
1{\sqrt{\kn}}\lr\frac{\bigl(\an\bigr)^{1/2}
|\qa|}{\si_{W_{(n)}}}\rr^3$, and by \eqref{1.7} it is a~magnitude
of the exact order $O\Bigl(\frac 1{\sqrt{\kn}}\Bigr)$, because
$\liminf_{\nty}\si^3_{W_{(n)}}>0$ under the condition $[A_1]$ when
$a_2<b_1$. Similarly for the third term we obtain the bound of the
order $O\Bigl(\frac 1{\sqrt{\mn}}\Bigr)$, whereas for the second term
we have
\begin{equation*}
\frac 1{\sqrt{n}\, \si^3_{W_{(n)}}}\,
\int_{\an}^{1-\bn}|F^{-1}(u)|^3\,du \le \frac
1{\sqrt{n}}\frac{|\,\qa\,|\vee
|\,\qb\,|}{\si_{W_{(n)}}}=O\Bigl(\frac 1{\sqrt{\kn}}+\frac
1{\sqrt{\mn}}\Bigr),
\end{equation*}
hence, for $\delta_{1,n}$ the desired estimate is valid. For
$\delta_{2,n}$ it follows directly from \eqref{qaqb}. For the
first term of $\delta_{3,n}$ by \eqref{qaqb} we obtain $\frac
1{\sqrt{n}}\bigl(\an\bigr)^{1/3}\lr\frac{\an}{f(\qa)\si_{W_{(n)}}}\rr^{5/3}\le
\frac
C{\sqrt{n}}\bigl(\an\bigr)^{1/3}\lr\bigl(\an\bigr)^{-1/2}\rr^{5/3}\le
\frac C{\sqrt{\kn}}$, where $C>0$ is some constant, independent
of~$n$ (cf~\eqref{qaqb}), and for the second term of
$\delta_{3,n}$ we similarly obtain the bound $O\Bigl(\frac
1{\sqrt{\mn}}\Bigr)$. Finally, for $\delta_{4,n}$ conditions
$[A'_2]$ and $[A'_3]$ directly yield $\frac
1{\sqrt{n}}\,\delta_{4,n}=O\Bigl(\frac 1{\sqrt{\kn}}+\frac
1{\sqrt{\mn}}\Bigr)$. The theorem is proved. $\quad \square$

\medskip
\noindent {\bf Proof of Theorem 1.3.} Also the validity of
this theorem is a simple consequence of Theorem~1.1. Take an arbitrary
$c>1/2$ and $A$, $B$ on the r.h.s. of \eqref{BE_1}, corresponding
to the value of $c$. Now to prove \eqref{BE_3} it suffices to repeat
the proof of Theorem~1.2, taking into account that $\an
\qa^2\vee\bn\qb^2\to 0$, as $\nty$ when $\si^2<\infty$. This gives
us the desired bound for $\frac
1{\sqrt{n}}\,(\delta_{1,n}+\delta_{2,n}+\delta_{3,n})$ at the
r.h.s. of \eqref{BE_1}. Finally, an application of the condition:
$\Psi_{\an,\frac 1{f(x)}}(B)=o\bigl((f(\qa)\ln\kn)^{-1}\bigr)$
and $\Psi_{1-\bn,\frac
1{f(x)}}(B)=o\bigl((f(\qb)\ln\mn)^{-1}\bigr)$ for every $B>0$, as
$\nty$, directly provides the bound $\frac
1{\sqrt{n}}\,\delta_{4,n}=o\Bigl(\frac 1{\sqrt{\kn}}+\frac
1{\sqrt{\mn}}\Bigr)$. The theorem is proved. $\quad \square$

\medskip
\noindent {\bf Proof~of~Corollary~1.1.} To prove this corollary we
apply \eqref{BE_1} with constants $A$ and $B$, corresponding to
$c=\frac 12$ in the general case, and to some value $c>1/2$ in
the special case where $\gamma =2 $ and $\si^2 <\infty$ to
obtain the bound \eqref{BE_3}. We must verify that the conditions $[A'_2]$
and $[A'_3]$ are satisfied in our case. Note that $F(x)$ and
$1-F(x)$ are regularly varying with index $\rho=-\gamma$
near $-\infty$ and $+\infty$ respectively. Moreover,
since  $\an=F(\qa)$ and $\bn=1-F(\qb)$, we find that
$\lim_{\nty}\frac{\an}{|\qa|f(\qa)}=lim_{x\to
-\infty}\frac{F(x)}{|x|f(x)}=\frac 1{\gamma}$ and
$\lim_{\nty}\frac{\bn}{|\qb|f(\qb)}=lim_{x\to
\infty}\frac{1-F(x)}{xf(x)}=\frac 1{\gamma}$ (cf.~Bingham et al.~\cite{bingham}), and the
conditions $[A''_2]$ (and hence, $[A'_2]$) is satisfied. This
implies that the quantity $\frac
1{\sqrt{n}}\bigl(\delta_{1,n}+\delta_{2,n}+\delta_{3,n}\bigr)$ has a~magnitude of the order $O\lr\frac 1{\sqrt{\kn}}+\frac
1{\sqrt{\mn}}\rr$. Moreover, in special case that $\gamma =2$ and
$\si^2<\infty$ the same quantity is of the smaller order, i.e.~$o\lr\frac 1{\sqrt{\kn}}+\frac
1{\sqrt{\mn}}\rr$ (cf.~proof of Theorems~1.2 and~1.3).

It remains to check that $[A'_3]$ holds true in our case. We will
verify that the~first inequality in $[A'_3]$ is satisfied,
for the second one we can apply similar argument. Set
$x_n=F^{-1}(\an)$ and $x_n+\triangle\,
x_n=F^{-1}\Bigl(\an+t\sqrt{\frac{\an\ln\kn}{n}}\Bigr)=
F^{-1}\Bigl(\an\bigl(1+t\sqrt{\frac{\ln\kn}{\kn}}\bigr)\Bigr)$,
where $|t|\le B$. Then $\Psi_{\an,\frac 1{f(x)}}\, (B)=\frac
1{f(x_n)}\mathfrak{D}_{\an,f}\, (B)$, where
$\mathfrak{D}_{\an,f}\, (B)=\sup_{|t|\le B}\lv \frac
{f(x_n)}{f(x_n+\triangle\, x_n)}-1\rv$, so that it remains to check
that $\mathfrak{D}_{\an,f}\, (B)=O\lr\frac 1{\ln\kn}\rr$ as
$\nty$.

Since the $df$ $F$ has an~unique inverse $F^{-1}$ on the set $U$, the
function $F^{-1}(\an)$ is regularly varying with index $-\frac
1{\gamma}$ when $\an\to 0$, and $|x_n|=(\an)^{-1/\gamma}L_1(\an)$,
where $L_1$ is a~slowly varying function when its argument tends to
zero. Therefore $\triangle\,
x_n=F^{-1}\Bigl(\an\bigl(1+t\sqrt{\frac{\ln\kn}{\kn}}\bigr)\Bigr)-
F^{-1}\Bigl(\an\Bigr)=|x_n|\lb
\lr1+t\sqrt{\frac{\ln\kn}{\kn}}\rr^{-1/\gamma} \frac
{L_1\bigl(\an\bigl(1+t\sqrt{\frac{\ln\kn}{\kn}}\bigr)\bigr)}{L_1\bigl(\an\bigr)}-1\rb$ with $L_1$ is as before and satisfying the requirement that it is in absolute
value of order $o(|x_n|)$. Then by the condition $[R]$ for every fixed
$t$ such that $|t|\le B$ we can write
\begin{equation}
\label{fxn} \frac {\lv f(x_n)-f(x_n+\triangle\,
x_n)\rv}{f(x_n+\triangle\, x_n)} =O\lr
\frac{\frac{f(x_n)}{|x_n|}|\triangle\, x_n|}{f(x_n)+O\lr
\frac{f(x_n)}{|x_n|}|\triangle\, x_n|\rr}  \rr =O \Bigl(
\lv\frac{\triangle\, x_n}{x_n}\rv\Bigr),
\end{equation}
as $\nty$. Next we note that $\triangle\,
x_n=\frac{1}{f\Bigl(F^{-1}\bigl( \an +\theta
t\an\sqrt{\frac{\ln\kn}{\kn}}\bigr)\Bigr)}t\an\sqrt{\frac{\ln\kn}{\kn}}$,
where  $0<\theta<1$, and by $[R]$ the latter quantity is equal
to $\frac{1}{f(x_n)+O\bigl(\frac{f(x_n)}{|x_n|}|\triangle\, x_n|
\bigr)}t\an\sqrt{\frac{\ln\kn}{\kn}}$. Then at the r.h.s. of
\eqref{fxn} we have a quantity of the order
$O\lr\frac{t\an\sqrt{\frac{\ln\kn}{\kn}}}{|x_n|f(x_n)(1+o(1))}\rr$,
as $\nty$. Since $|x_n|f(x_n)\sim F(x_n)\frac 1{\gamma}$ due to
the regularly varying property, and because
$F(x_n)=F(F^{-1}(\an))=\an$, we obtain that the quantity at the
r.h.s. of \eqref{fxn} is of the order $O\lr
\sqrt{\frac{\ln\kn}{\kn}}\rr$ uniformly in all $|t|\le B$. This
implies that $\Psi_{\an,\frac
1{f(x)}}(B)=o\bigl(\bigl(f(\qa)\ln\kn\bigr)^{-1}\bigr)$, and
similarly we obtain that $\Psi_{1-\bn,\frac
1{f(x)}}(B)=o\bigl(\bigl(f(\qb)\ln\mn\bigr)^{-1}\bigr)$. We can conclude that
under
our conditions  
we have $\frac
1{\sqrt{n}}\bigl(\delta_{4,n}\bigr)=o\lr\frac 1{\sqrt{\kn}}+\frac
1{\sqrt{\mn}}\rr$, what completes the proof of the corollary.
 $\quad \square$

\medskip
\noindent {\bf Proof~of~Corollary~1.2.}  It suffices to check that the
conditions of Theorem~1.2 are satisfied. Condition $[A_1]$, inequalities $0<a_1<b_2<1$
and continuity of $f$ together imply that both $\liminf_{\nty}\bigl(f(\qa)\wedge
f(\qb)\bigr)>0$ and $\liminf_{\nty}\si_{W_{(n)}}>0$ are automatically satisfied. Hence,
$[A'_2]$ holds true. Moreover, by the H\"{o}lder condition of
degree $d>0$ we have: $\Psi_{\an,\frac 1{f(x)}}(B)=O\lr \frac
{\an\ln \kn}{n}\rr^{d/2}=O\lr\frac {\ln n}{n}\rr^{d/2}= o(\frac 1{\ln n})$, as $\nty$. The same argument
is valid for $\Psi_{1-\bn,\frac 1{f(x)}}(B)$. Thus, $[A'_3]$ also holds true.  The
corollary is proved. $\quad \square$

\medskip
\noindent {\bf Proof~of~Theorem~1.4.} Similarly as in the proof of Theorem 1.1
(cf.~\eqref{petr_l}) we write
\begin{equation}
\label{petr_Tm5}
\begin{split}
 F_{U,n}\Bigl(x-\frac{\Delta_n}{\si_{W_{(n)}}}\Bigr)-\boldsymbol{P}\lr |R_n|> \Delta_n \rr&\le
F_{T_n}\Bigl(x+\frac{b_n}{\si_{W_{(n)}}}\Bigr)\\ &\le
F_{U,n}\Bigl(x+\frac{\Delta_n}{\si_{W_{(n)}}}\Bigr)+\boldsymbol{P}\lr |R_n|>
\Delta_n \rr,
\end{split}
\end{equation}
where $\Delta_n$ is as in Lemma~2.1, $b_n$ is as in \eqref{G_n},
and $F_{U,n}$ is the $df$ of $\bigl(L_n+U_n\bigr)/\si_{W_{(n)}}$
(cf.~\eqref{petr_l}) --- the $U$-statistic of degree~2 with
a~kernel which depends on $n$. Our smoothness condition $[A_1]$
implies that the $df$ of the r.v. $W_1(n)$ has a~positive density
on a~Borel set $D_n\subset \mathbb{R}$ with nonzero Lebesgue
measure and that there exists an integer $n_0$ such that
$\underset{n=n_0}{\overset{\infty}{\bigcap}} D_n\supset D\neq
\emptyset$, and the set $D$ has nonzero Lebesgue measure. Hence,
the members of the consequence of the $df$'s of r.v.'s $W_1(n)$
have a~common nontrivial absolutely continuous component
independent of $n$ for all sufficiently large $n$. This yields
\begin{equation}
\label{Cramer_Tm5}
\limsup_{\nty}\limsup_{|t|\to\infty}\bigl|\boldsymbol{E}\exp\bigl(i\,
t\sqrt{n}\, L_{n,1}\bigr)
\bigr|=\limsup_{\nty}\limsup_{|t|\to\infty}\bigl|\boldsymbol{E}\exp\bigl(i\,
t\, W_1(n)\bigr) \bigr| <1,
\end{equation}
hence the sequence of the first canonical functions of the
$U$-statistic satisfies the Cramer condition, and we can apply
a~result by Bentkus et al.~\cite{bgz}. The one term Edgeworth expansion of the $df$
$F_{U,n}(x)=\boldsymbol{P}\Bigl(\bigl(L_n+U_n\bigr)/\si_{W_{(n)}}\le x\Bigr)$
is equal to $G_{U,n}(x)=\Phi(x)-\frac{\phi(x)}{6\sqrt{n}}\bigl(
\lambda_{1_{(n)}}+3\lambda_{2_{(n)}}\bigr)(x^2-1)$ (cf.~Section~2,
cf.~also Bentkus et al.~\cite{bgz}, page~855). Write
\begin{equation}
\label{Tm5_14} \sup_{x\in\mathbb{R}}\Bigl| F_{U,n}\Bigl(
x\pm\frac{\Delta_n}{\si_{W_{(n)}}}\Bigl)- G_{U,n}(x)\Bigr|\le
\Delta_{n,1}+\Delta_{n,2},
\end{equation}
where $\Delta_{n,1}=\sup_{x\in\mathbb{R}}\bigl| F_{U,n}\bigl(
x\bigr) - G_{U,n}(x)\bigr|$, \
$\Delta_{n,2}=\sup_{x\in\mathbb{R}}\Bigl| G_{U,n}\Bigl(
x\pm\frac{\Delta_n}{\si_{W_{(n)}}}\Bigl)- G_{U,n}(x)\Bigr|$. To
estimate $\Delta_{n,1}$ we apply Theorem~1.2 of Bentkus et
al.~\cite{bgz}, taking into account the Remark~1.3, given on page~856
in cited paper. Then we obtain
\begin{equation}
\label{Tm5_15} \Delta_{n,1}\le \frac Cn \, \Bigl(\,
\frac{\boldsymbol{E}(W_1(n))^4}{\si^4_{W_{(n)}}}+\frac{\gamma_{2+\varepsilon}}
{\si^{2+\varepsilon}_{W_{(n)}}} \Bigr),
\end{equation}
where $\varepsilon>0$ is an arbitrary constant, the constant $C>0$
depends on $\varepsilon$ and does not depend on $n$ (Note that the
quantity $\Delta^2_3$ appearing in Theorem~1.2 is zero in case of a~$U$-statistic of degree~2 (cf.~Bentkus et al.~\cite{bgz}, page~858)),
and
\begin{equation}
\label{Tm5_16} \begin{split} \gamma_{2+\varepsilon} &=\boldsymbol{E}\lv
n^{3/2}U_{n,(1,2)}\rv^{2+\varepsilon}\\
&\le 2^{1+\varepsilon}\lb\frac{\bigl(\boldsymbol{E}\,
|I_{\an}(X_1)-\an|^{2+\varepsilon}\bigr)^2}
{f^{2+\varepsilon}(\qa)}+ \frac{\bigl(\boldsymbol{E}\,
|I_{1-\bn}(X_1)-(1-\bn)|^{2+\varepsilon}\bigr)^2}
{f^{2+\varepsilon}(\qb)}\rb\\
&<2^{1+\varepsilon}\lb\alp_n^2\frac{1} {f^{2+\varepsilon}(\qa)}+
\be_n^2\frac{1} {f^{2+\varepsilon}(\qb)}\rb.
\end{split}
\end{equation}
Relations \eqref{Tm5_15} - \eqref{Tm5_16} imply that
$\Delta_{n,1}\le \frac {C_1}n
\bigl(\delta_{1,n}+\delta_{2,n}\bigr)$, and since $G'_{U,n}(x)$ is
bounded uniformly in \ $x$, we obtain $\Delta_{n,2}\le  {C_2} \,
\frac {\Delta_n}{\si_{W_{(n)}}}=\frac {C_2}{n^{3/4}} \delta_{4,n}+
\frac {C_2}{n^{1/2}}\delta_{5,n}$, where $C_i$, $i=1,2$, some
positive constants not depending on $n$. These estimates,
relation \eqref{petr_Tm5} and Lemma~2.1 together imply that
\begin{equation}
\label{Tm5_17}  \sup_{x\in\mathbb{R}}\Bigl| F_{T_n}(x) -
G_{U,n}\Bigl(x-\frac{b_n}{\si_{W_{(n)}}}\Bigr)\Bigr|\le \frac
{C_1}n \bigl(\delta_{1,n}+\delta_{2,n}\bigr)+\frac {C_2}{n^{3/4}}
\delta_{4,n}+ \frac
{C_2}{n^{1/2}}\delta_{5,n}+C_3\lr\kn^{-c}+\mn^{-c}\rr.
\end{equation}
It remains to note that since $G'_{U,n}(x)$ and $G''_{U,n}(x)$ are
bounded uniformly in $x$, we have
\begin{equation}
\label{Tm5_18}  \sup_{x\in\mathbb{R}}\Bigl|
G_{U,n}\Bigl(x-\frac{b_n}{\si_{W_{(n)}}}\Bigr)-G_{U,n}(x)
\Bigr|\le C \Bigl(\frac{|
(\lambda_{1_{(n)}}+3\lambda_{2_{(n)}})\,b_n|}{\sqrt{n}\,\si_{W_{(n)}}}
+\frac{b^2_n}{\si^2_{W_{(n)}}}\Bigr) \le C\Bigl( \frac
{\delta_{6,n}}{n^{1/2}}+\frac {\delta_{3,n}}{n}\Bigr),
\end{equation}
where $C$ is some constant not depending on $n$ . Relations
\eqref{Tm5_17} -- \eqref{Tm5_18} imply \eqref{EE_n1}. The theorem
is proved. $\quad \square$

\medskip
\noindent {\bf Proof~of~Corollary~1.3.} To prove this corollary we
apply Theorem~1.4 with $c=3/4$, and check that the quantity at the
r.h.s. of \eqref{EE_n1} is of the order $O\Bigl(\frac
{(\ln\kn)^{5/4}}{\kn^{3/4}}+\frac
{(\ln\mn)^{5/4}}{\mn^{3/4}}\Bigr)$ under conditions of
Corollary~1.3. Similarly as in proof of Theorem~1.2 using condition
$[A'_2]$ and \eqref{1.7} we easily verify that $\frac 1n
\bigl(\delta_{1,n}+\delta_{3,n}\bigr)=O\Bigl(\frac 1{\kn}+\frac
1{\mn}\Bigr)$. By condition $[A'_2]$ for $\frac 1n\delta_{2,n}$ we
have the bound $O\Bigl(\frac
1n\Bigl[\an^{-1-3\varepsilon/2}+\bn^{-1-3\varepsilon/2}\Bigr]\Bigr)=
O\Bigl(\frac 1{\kn}\Bigl(\frac
{n}{\kn}\Bigr)^{3\varepsilon/2}+\frac 1{\mn}\Bigl(\frac
{n}{\mn}\Bigr)^{3\varepsilon/2}\Bigr)$, and by condition $[L]$
the latter quantity is of the order $O\Bigl(\frac
{1}{\kn^{3/4}}+\frac {1}{\mn^{3/4}}\Bigr)$ if $s\ge
6\varepsilon/(1+6\varepsilon)$. For $\frac 1{n^{3/4}}\delta_{2,n}$
we have the desired bound directly by $[A'_2]$, for $\frac
1{n^{1/2}}\delta_{5,n}$ we get the same bound directly by the
conditions $[A'_2]$ and the condition on $\Psi_{\an, \frac
1{f(x)}}(B)$ and  $\Psi_{1-\bn, \frac 1{f(x)}}(B)$. To treat $\frac
1{n^{1/2}}\delta_{6,n}$ we use the same argument as before based on
\eqref{1.7} and condition $[A'_2]$, which leads to a~bound of
order $O\Bigl(\frac 1{\kn}+\frac 1{\mn}\Bigr)$ for this term. The
corollary is proved. $\quad \square$

\medskip
\noindent {\bf Proof~of~Theorem~1.5.} Similarly as in proof of
Corollary~1.1 we check that the condition $[A''_2]$ (and hence
$[A'_2]$) is satisfied. Moreover, conditions for $\Psi_{\an, \frac
1{f(x)}}(B)$  and  $\Psi_{1-\bn, \frac 1{f(x)}}(B)$ (cf.~ corollary~1.3) are  satisfied if the
conditions $[A_1]$, $[R]$ hold true (cf.~proof of the
Corollary~1.1). Thus, we obtain the validity of \eqref{EE_n2} as
a~consequence of Corollary~1.3.  The theorem is proved. $\quad \square$

\medskip
\noindent {\bf Proof~of~Corollary~1.4.} This corollary follows
directly from \eqref{EE_n1}. Indeed, in our conditions we have
$\frac
1n\bigl(\delta_{1,n}+\delta_{2,n}+\delta_{3,n}\bigr)=O\lr\frac
1n\rr$, \ $n^{-3/4}\delta_{4,n}=O\lr\frac{\ln^{5/4}
n}{n^{3/4}}\rr$, \ $n^{-1/2}\delta_{6,n}=O\lr\frac 1n\rr$, and by
H\"{o}lder condition  $\ln n\bigl(\Psi_{\an,\frac
1{f(x)}}(B)+\Psi_{1-\bn,\frac 1{f(x)}}(B)\bigr)=O\lr\ln n\bigl(
\frac{\ln n}{n}\bigr)^{d/2}\rr=o\bigl(n^{-d/2+\varepsilon}\bigr)$
for every $\varepsilon>0$. These bounds imply that
$n^{-3/4}\delta_{4,n}+n^{-1/2}\delta_{5,n}
=o\bigl(n^{-1/2-p}\bigr)$ for every $p<\min(1/4,d/2)$. The
corollary is proved. $\quad \square$

\medskip
\noindent {\bf Proof~of~Theorem~1.6.} First we write
$F_{T_n,S}(x)=\boldsymbol{P}\Bigl(\frac{L_n+U_n+b_n}{\widehat{\si}_{W_{(n)}}}+
\frac{R_{n_{\phantom{n_n}}}}{\si_{W_{(n)}}}\,
\frac{\si_{W_{(n)}}}{\widehat{\si}^{\phantom{n^n}}_{W_{(n)}}}\le
x\Bigr)$, where by Lemma~2.1: \
$p_{n,1}:=\boldsymbol{P}\bigl(|R_n|>\Delta_n\bigr)=O\bigl(\kn^{-c}+\mn^{-c}\bigr)$,
for every $c>0$, and $\Delta_n$ is as in \eqref{L4.1_1}. By
Lemma~2.2 the main term of the quantity
$\frac{\widehat{\si}^2_{W_{(n)}}}{\si^2_{W_{(n)}}}-1$ \ is \
$\frac{V_n}{\si^2_{W_{(n)}}}$, for which by Chebyshev's inequality
for every $t>0$ we have
$p_{n,2}=\boldsymbol{P}\bigl(\frac{|V_n|}{\si^2_{W_{(n)}}}
>2t\bigr)\le \frac{\boldsymbol{E}V^2_n}{4t^2 \si^4_{W_{(n)}}} \le
\frac C{t^2}\bigl(\frac{\boldsymbol{E}(W_1(n))^4}{n\, \si^4_{W_{(n)}}}
+q^2_{\an}+q^2_{\bn}\bigr)$, where $C>0$ is some constant
independent of $n$ and $t$ (cf.~\eqref{EV_n}), and because
\eqref{1.7} the latter quantity is of the order $O\bigl(\frac
1{\kn}+\frac
1{\mn}+q^2_{\an}+q^2_{\bn}\bigr)=o(\delta_{1,S}(n)+\delta_2(n))$,
where $\delta_{1,S}(n)$, $\delta_2(n)$ as in~\eqref{EE_s1}
and~\eqref{delta_s} respectively. This implies that $\Bigl|
\frac{\si_{W_{(n)}}}{\widehat{\si}^{\phantom{n^n}}_{W_{(n)}}}-1\Bigr|\le
t$ with probability of the order $p_{n,2}$. \ Put
$P_n=p_{n,1}+p_{n,2}$. Then we obtain
\begin{equation}
\label{Th8_p1}
\widetilde{F}_{U_n,S}\Bigl(x-\frac{\Delta_n(1+t)}{\si_{W_{(n)}}}\Bigr)-P_n
\, \le \, F_{T_n,S}(x) \, \le \,
\widetilde{F}_{U_n,S}\Bigl(x+\frac{\Delta_n(1+t)}{\si_{W_{(n)}}}\Bigr)+P_n,
\end{equation}
where
\begin{equation}
\label{F_s_tilde}
\widetilde{F}_{U_n,S}(x)=\boldsymbol{P}\Bigl(\frac{L_n+U_n+b_n}{\si_{W_{(n)}}}\le
x\Bigl(1+\frac{V_n}{\si^2_{W_{(n)}}}+R_{n,S}\Bigr)^{1/2}\Bigr),
\end{equation}
where $R_{n,S}$ is the remainder term from Lemma~2.2. Note that
$1+\frac{V_n}{\si^2_{W_{(n)}}}>0$ for all sufficiently large $n$
with probability of the order $P_n$, and
$|R_{n,S}|=O(\Delta_{n,S})$ with probability
$1-O(\kn^{-c}+\mn^{-c})$, where $\Delta_{n,S}$ is as in
\eqref{delta_s}.

Since $H'_n(x)\, x$ is bounded from above uniformly in $x$, it is
enough to prove that $H_n(x)$ is the expansion for the r.h.s of
\eqref{F_s_tilde} without $R_{n,S}$, because omitting of it gives
a~remainder term of the order $O(\Delta_{n,S})$, which presences
at the r.h.s. of \eqref{EE_s1}. Write
\begin{equation}
\label{Th8_p2}
\begin{split}
&\boldsymbol{P}\Bigl(\frac{L_n+U_n+b_n}{\si_{W_{(n)}}}\le
x\Bigl(1+\frac{V_n}{\si^2_{W_{(n)}}}\Bigr)^{1/2}\Bigr)\\
=&\boldsymbol{P}\Bigl(\frac{L_n+U_n+b_n}{\si_{W_{(n)}}}- x\Bigl{\{}
\Bigl(1+\frac{V_n}{\si^2_{W_{(n)}}}\Bigr)^{1/2}-1\Bigr\}\le
x\Bigr).
\end{split}
\end{equation}
Since $\frac{|V_n|}{\si^2_{W_{(n)}}}\le t$ with probability
$1-P_{n,1}$ for every $t>0$, we can apply as in Putter and
van~Zwet~\cite{pz} (cf.~also~\cite{gh2006}-\cite{gh2007}) the
following inequality: $1+\frac z2-\frac {z^2}4\le (1+z)^{1/2}\le
1+ \frac z2$, \ $|z|\le 4/5$. Note that
$\frac{V^2_n}{4\si^4_{W_{(n)}}}\le \frac 12 \Bigl(
\frac{V^2_{n,1}}{\si^4_{W_{(n)}}}+\frac{V^2_{n,2}}{\si^4_{W_{(n)}}}\Bigr)$,
where $V^2_{n,i}$, $i=1,2$, are as in \eqref{V_n}, and note that
$\frac{V^2_{n,1}}{\si^4_{W_{(n)}}}=O(\ln\kn q^2_{\an}+\ln\mn
q^2_{\bn})$ (cf.~\eqref{EV_n}), the latter quantity contributes to
 to $\delta_{n,S}$  on the r.h.s. of
\eqref{EE_s1} (because it is a term of $\Delta_{n,S}$). It follows
that we have to show that
\begin{equation}
\sup_{x\in\mathbb{R}}\Bigl|\boldsymbol{P}\Bigl(\frac{L_n+U_n+b_n}{\si_{W_{(n)}}}-\frac{x\,V_n}
{2\si^2_{W_{(n)}}}\le x\Bigr)-H_n(x)\Bigr|\le
C(\delta_n+\delta_{n,S}),
\end{equation}
\begin{equation}
\sup_{x\in\mathbb{R}}\Bigl|\boldsymbol{P}\Bigl(\frac{L_n+U_n+b_n}{\si_{W_{(n)}}}-\frac{x\,V_n}
{2\si^2_{W_{(n)}}}+\frac{x\,V^2_{n,2}} {2\si^4_{W_{(n)}}}\le
x\Bigr)-H_n(x)\Bigr|\le C(\delta_n+\delta_{n,S}),
\end{equation}
where $C>0$ is some constant independent of $n$ and $\delta_n$ is
as in \eqref{EE_s1}. Define
$\widetilde{H}_n(x)=H_n(x)+\frac{b_n\phi(x)}{\si_{W_{(n)}}}$
(i.e.~$\widetilde{H}_n(x)$ is $H_n(x)$ without bias term). Note
that $\frac{|b_n|}{\si_{W_{(n)}}}\le \frac 12(q_{\an}+q_{\bn})$.
Since $x\,H'_n(x)$ and $x^2\,H'_n(x)$ are bounded, we obtain:
$H_n(x+\frac{b_n}{\si_{W_{(n)}}})=H_n(x)+\phi(x)\frac{b_n}{\si_{W_{(n)}}}
+O\Bigl(\frac{|\lambda_{1_{(n)}}|+|\lambda_{2_{(n)}}|}{\sqrt{n}}\,
\frac{|b_n|}{\si_{W_{(n)}}}+b^2_n\Bigr)=H_n(x)+\phi(x)\frac{b_n}{\si_{W_{(n)}}}
+O(\frac 1{\sqrt{n}}\delta_{6,n}+q^2_{\an}+q^2_{\bn})=
H_n(x)+\phi(x)\frac{b_n}{\si_{W_{(n)}}}+O(\delta_n+\delta_{n,S})$.
It follows that we should prove that
\begin{equation}
\label{Th8_p3}
\sup_{x\in\mathbb{R}}\Bigl|\boldsymbol{P}\Bigl(\frac{L_n+U_n}{\si_{W_{(n)}}}-\frac{x\,V_n}
{2\si^2_{W_{(n)}}}\le x\Bigr)-\widetilde{H}_n(x)\Bigr|\le C(
\delta_n+\delta_{n,S}),
\end{equation}
\begin{equation}
\label{Th8_p4}
\sup_{x\in\mathbb{R}}\Bigl|\boldsymbol{P}\Bigl(\frac{L_n+U_n}{\si_{W_{(n)}}}-\frac{x\,V_n}
{2\si^2_{W_{(n)}}}+\frac{x\,V^2_{n,2}} {2\si^4_{W_{(n)}}}\le
x\Bigr)-\widetilde{H}_n(x)\Bigr|\le C(\delta_n+\delta_{n,S}).
\end{equation}

First we prove \eqref{Th8_p3}. Since $V_n$ is a sum of centered
i.i.d. r.v.'s, we obtain that
$U_x=\frac{L_n+U_n}{\si_{W_{(n)}}}-\frac{x\,V_n}
{2\si^2_{W_{(n)}}}$ is a~centered $U$-statistic of degree two, and
as in proof of Theorem~1.4 we find that in view of our smoothness
assumption $[A_1]$ the  Cramer condition is satisfied. Put
$\nu_n=\ln(\kn\wedge\mn)$. First we prove that \eqref{Th8_p3}
holds true uniformly in $x$: $|x|<\nu_n$. By Theorem~1.1 of
Bentkus et al.~\cite{bgz} (taking into account the Remark~1.3 given on
page~856 in cited paper) after simple computation of the fourth
moment of $U_x$ we obtain
\begin{equation}
\label{Th8_p5} \begin{split}
\sup_{|x|<\nu_n}\Bigl|\boldsymbol{P}\Bigl(U_x\le
x\Bigr)-\widetilde{G}_n(x)\Bigr|&\le \frac
Cn\lr\frac{\boldsymbol{E}(W_1(n))^4}{\si^4_{W_{(n)}}}+\frac{\nu_n^4}{n^2}\lb
\frac{\boldsymbol{E}(W_1(n))^8}{\si^8_{W_{(n)}}}\rp\rp\\
&\lp\lp+\frac{\alp_n^5\qa^4}{f^4(\qa)\si^8_{W_{(n)}}}+
\frac{\be_n^5\qb^4}{f^4(\qb)\si^8_{W_{(n)}}}\rb +
\frac{\gamma_{2+\varepsilon}}{\si^{2+\varepsilon}_{W_{(n)}} }\rr,
\end{split}
\end{equation}
where $\varepsilon>0$ is an arbitrary independent of $n$ constant
(which depends on $\varepsilon>0$, cf.~Bentkus et al.~\cite{bgz},
page~856), $\gamma_{2+\varepsilon}$ is as in proof of Theorem~1.4,
and
\begin{equation}
\label{Th8_p6} \widetilde{G}_n(x)=\Phi\Bigl(\frac{x}{\si_x}\Bigr)
-\frac{k_{3,x}}{6\si^3_x}\lb\Bigl(\frac{x}{\si_x}\Bigr)^2
-1\rb\phi\Bigl(\frac{x}{\si_x}\Bigr),
\end{equation}
where $\si^2_x=\boldsymbol{E}U^2_x$, \ $k_{3,x}=\boldsymbol{E}U^3_x$.
Relation \eqref{1.7} implies that $\frac{\nu_n^4}{n^3}\lb
\frac{\boldsymbol{E}(W_1(n))^8}{\si^8_{W_{(n)}}}\rp$ $\lp
+\frac{\alp_n^5\qa^4}{f^4(\qa)\si^8_{W_{(n)}}}+
\frac{\be_n^5\qb^4}{f^4(\qb)\si^8_{W_{(n)}}}\rb=O\lr\nu_n^4(\frac
1{\kn^3}+\frac 1{\mn^3})+ \frac{\nu_n^4}{n} \lb\frac
1{\an}q_{\an}^4+\frac
1{\bn}q_{\bn}^4\rb\rr=o(\delta_{1,S}(n)+\delta_{2,S}(n))$, where
$\delta_{i,S}(n)$ is as in \eqref{EE_s1}. Moreover, as in proof of
Theorem~1.4 we obtain that
$\frac{\gamma_{2+\varepsilon}}{\si^{2+\varepsilon}_{W_{(n)}}
}=O\lr\delta_{2,n}\rr$. Thus, at the r.h.s. of \eqref{Th8_p5} we have
desired bound  $C\lr \frac 1n\lb\delta_{1,n}+
\delta_{2,n}\rb +\delta_{1,S}(n)\rr\le C(\delta_n+\delta_{n,S})$
(here $\delta_{i,n}$, $i=1,2$, are two terms of $\delta_n$,
cf.~\eqref{EE_n1}).

Next consider $\widetilde{G}_n(x)$. We have
$\si^2_x=\boldsymbol{E}U^2_x=\boldsymbol{E}\lr\frac{L_n+U_n}{\si_{W_{(n)}}}\rr^2-\frac
x{\si^2_{W_{(n)}}}\boldsymbol{E}\bigl(\,[L_n+U_n]\,V_n\bigr)+\frac{x^2}{4}\boldsymbol{E}\lr
\frac{V^2_n}{\si^4_{W_{(n)}}}\rr$. Since $U_n$ and $V_n$ are
uncorrelated, after simple computations using formulas
\eqref{4.1}-\eqref{4.5} and \eqref{V_n} we obtain
$\boldsymbol{E}\bigl([L_n+U_n]V_n\bigr)=\boldsymbol{E}\bigl(L_n
V_n\bigr)=\frac
1{\sqrt{n}}\bigl(\gamma_{3,W_{(n)}}+2\delta_{2,W_{(n)}}\bigr)$,
and hence (cf.~\eqref{4.6},~\eqref{EV_n}),
\begin{equation}
\label{si_x}
\si^2_x=1-\frac{x(\lambda_{1_{(n)}}+2\lambda_{2_{(n)}})}{\sqrt{n}}
+O\Bigl(\nu_n^2\Bigl[\frac{\boldsymbol{E}(W_1(n))^4}{n\si^4_{W_{(n)}}}
+q^2_{\an}+q^2_{\bn}\Bigr]\Bigr).
\end{equation}
Moreover, relations \eqref{4.1}-\eqref{4.5}, \eqref{4.9},
\eqref{V_n} and \eqref{1.7} after simple computations yield
\begin{equation}
\label{k3_x} k_{3,x}=\boldsymbol{E}
U_x^3=\frac{\lambda_{1_{(n)}}+3\lambda_{2_{(n)}}}{\sqrt{n}}
+O(\delta_{n,S}).
\end{equation}
Note that estimating of the remainder term at the r.h.s. of
\eqref{k3_x} is  essentially based on relation \eqref{1.7}, which
we use to bound the moments $\boldsymbol{E}W^r_1(n)/\si^r_{W_{(n)}}$
(cf.~proof of Theorem~1.2), where the largest power appearing here
is $r=6$. The relations \eqref{Th8_p6} -- \eqref{k3_x} together
imply that
\begin{equation}
\label{Th8_p7} \widetilde{G}_n(x)=\Phi\Bigl(\frac{x}{\si_x}\Bigr)
-\frac{\lambda_{1_{(n)}}+3\lambda_{2_{(n)}}}{6\,\sqrt{n}}
\bigl(x^2-1\bigr)\phi\bigl(x \bigr)+O(\delta_{n,S}),
\end{equation}
for $|x|\le \nu_n$, that is $\si_x$ influences the first term
of the expansion only through the term $\Phi\Bigl(\frac{x}{\si_x}\Bigr)$
(cf.~proof of Theorem~1.2 in Putter \& van~Zwet~\cite{pz}). Then
using \eqref{si_x} and \eqref{1.7}, we obtain:
$\Phi\Bigl(\frac{x}{\si_x}\Bigr)=\Phi(x)+\phi(x)\frac{x^2(\lambda_{1_{(n)}}+
2\lambda_{2_{(n)}})}{2\,\sqrt{n}} +O(\delta_{n,S})$. Thus,
$\widetilde{G}_n(x)=\widetilde{H}_n(x)+O(\delta_{n,S})$ for
$|x|\le \nu_n$. Then we argue as Putter \&
van~Zwet~\cite{pz} (cf.~also~\cite{gh2006}-\cite{gh2007}): for
$x<-\nu_n$ \ $\widetilde{H}_n(x)=O(\kn^{-c}+\mn^{-c})$, and for
$x>\nu_n$ \ $1-\widetilde{H}_n(x)=O(\kn^{-c}+\mn^{-c})$, where
$c>0$ is an arbitrary constant. So, monotonicity of a distribution
function implies \eqref{Th8_p3}.

It remains to prove \eqref{Th8_p4}. As before we see that it is
enough to prove it taking supremum in $x:\ |x|<\nu_n$. We must
prove that the presence of $\frac{xV^2_{n,2}}{2\si^4_{W_{(n)}}}$
does not influence on the expansion and the order of the bound at
the r.h.s. of \eqref{Th8_p4}. Note that $xV^2_{n,2}=\frac
x{n^2}\Bigl(\sum_{i=1}^{n}\bigl[(W_i(n)-\mu_{W_{(n)}})^2-
\si^2_{W_{(n)}}\bigr]\Bigr)^2$ is a $U$-statistic of degree two,
and its Hoeffding's decomposition is
\begin{equation*}
\begin{split}
&\frac{x}{n}\boldsymbol{E}\bigl[(W_1(n)-\mu_{W_{(n)}})^2-
\si^2_{W_{(n)}}\bigr]^2+\frac{x}{n^2}\sum_{i=1}^{n}\lr
\bigl[(W_i(n)-\mu_{W_{(n)}})^2- \si^2_{W_{(n)}}\bigr]^2
-\boldsymbol{E}\bigl[(W_1(n)-\rp\\
&\lp \mu_{W_{(n)}})^2 -
\si^2_{W_{(n)}}\bigr]^2\rr+\frac{2x}{n^2}\sum_{1\le i<}\sum_{j\le
n}\bigl[(W_i(n)-\mu_{W_{(n)}})^2-
\si^2_{W_{(n)}}\bigr]\bigl[(W_j(n)-\mu_{W_{(n)}})^2-
\si^2_{W_{(n)}}\bigr].
\end{split}
\end{equation*}
Since
$\frac{x}{n\si^4_{W_{(n)}}}\boldsymbol{E}\bigl[(W_1(n)-\mu_{W_{(n)}})^2-
\si^2_{W_{(n)}}\bigr]^2=O\lr\frac{\nu_n}{n}\frac{\boldsymbol{E}W^4_1(n)}
{\si^4_{W_{(n)}}}\rr$ and as $\widetilde{H}'(x)$ is bounded
uniformly in $x\in\mathbb{R}$, the constant term of Hoeffding's
decomposition contributes to a remainder term and can be omitted.
Since $\widetilde{U}_x=\frac{L_n+U_n}{\si_{W_{(n)}}}-\frac{x\,V_n}
{2\si^2_{W_{(n)}}}+\frac{x\,V^2_{n,2}}
{2\si^4_{W_{(n)}}}-\frac{x}{2n\si^4_{W_{(n)}}}\boldsymbol{E}\bigl[(W_1(n)-\mu_{W_{(n)}})^2-
\si^2_{W_{(n)}}\bigr]^2$ is a centered $U$-statistic of degree
two, and as $\frac{x\,V^2_{n,2}} {2\si^4_{W_{(n)}}}$ is even less
than $\frac{x\,V_n} {2\si^2_{W_{(n)}}}$ for all sufficiently large
$n$, applying Theorem~1.1 of Bentkus et al.~\cite{bgz} we come to the
same estimate as in \eqref{Th8_p5} with
$\widetilde{G}_n(x)=\Phi\Bigl(\frac{x}{\
\widetilde{\si}^{\phantom{a^b}}_x}\Bigr)
-\frac{\lambda_{1_{(n)}}+3\lambda_{2_{(n)}}}{6\,\sqrt{n}}
\bigl(x^2-1\bigr)\phi\bigl(x \bigr)+O(\delta_{n,S})$, where for
$\widetilde{\si}^2_x=\boldsymbol{E}\widetilde{U}^2_x$ after some
simple but rather tedious computations we obtain:
$\widetilde{\si}^2_x=1-\frac{x(\lambda_{1_{(n)}}+
2\lambda_{2_{(n)}})}{\sqrt{n}}
+O\Bigl(\nu_n^2\Bigl[\frac{\boldsymbol{E}(W_1(n))^4}{n\si^4_{W_{(n)}}}
+q^2_{\an}+q^2_{\bn}\Bigr]\Bigr)+o(\delta_{n,S})$. Thus, as well
as before $\Phi\Bigl(\frac{x}{\
\widetilde{\si}^{\phantom{a^b}}_x}\Bigr)=\Phi(x)+\phi(x)\frac{x^2(\lambda_{1_{(n)}}+
2\lambda_{2_{(n)}})}{2\,\sqrt{n}} +O(\delta_{n,S})$, and
$\widetilde{G}_n(x)=\widetilde{H}_n(x)+O(\delta_{n,S})$. This
implies that \eqref{Th8_p4} is valid when supremum is taken all
over $x:\ |x|<\nu_n$, then as before by a~monotonicity argument
we obtain the validity of it on the whole real line. The theorem
is proved.$\quad \square$

\medskip
\noindent {\bf Proof~of~Corollary~1.5.} To prove this corollary we
apply Theorem~1.6 with $c=3/4$, and check that
$\delta_n+\delta_{n,S}=O\Bigl(\frac
{(\ln\kn)^{5/4}}{\kn^{3/4}}+\frac
{(\ln\mn)^{5/4}}{\mn^{3/4}}\Bigr)$, as $\nty$. We obtain the desired bound for $\delta_n$ similarly
as in the proof of Corollary~1.3. The check for
$\delta_{n,S}$ is similar. The corollary is proved.$\quad \square$

\medskip
\noindent {\bf Proof~of~Theorem~1.7.} The proof is similar to the proof of Theorem~1.5. $\quad \square$

\section{Some Bahadur -- Kiefer type representations}
\label{Bahadur}
In this section we state and prove two lemmas used in
our proofs. In essence this lemmas extend corresponding auxiliary results obtained in \cite{gh2006}, \cite{gh2007}
for a~special case of the central sample  quantiles  $\xi_{\alp n:n}$ ($0<\alp<1$ is fixed)  to the case that $\an$ is a~sequence, which in particular
 can tend  to $0$ or to $1$ (i.e. to the case of intermediate sample  quantiles).

Let $\kn$ be a sequence of positive integers such that $\kn \to
\infty$, recall that $\an=\kn/n$, $0\le \liminf \an\le\limsup\an<1$, $\xi_{\alp_n n:n}=F_n^{-1}(\alp_n)$ denote the corresponding sample  quantiles,
and let $U_a$ be the set defined in \eqref{1.6}.
Let $G(x)$, $x\in R$, be a real-valued function, $g=G'$ -- its derivative when it exists,
and let $(g/f)(x)$ denote the ratio $g(x)/f(x)$, $(|g|/f)(x)$ --- the ratio $|g(x)|/f(x)$.
\begin{lemma}
\label{lem_A.1}
Suppose that $F^{-1}$ and $G$ are
differentiable on the sets $U_a$ and $F^{-1}(U_a)$ respectively.
Then
\begin{equation}
\label{(A.1)}
G(\qan)-G(\qa)=-[F_n(\qa)-F(\qa)]\frac gf(\qa)+R_n,
\end{equation}
where  $\boldsymbol{P}(|R_n|>\Delta_n)=O\lr\kn^{-c}\rr$ for each $c>0$, and
\begin{equation*}
\Delta_n=A\, \an\lb \frac
{|g|}f(\qa)\lr\frac{\ln\kn}{\kn}\rr^{3/4}+\Psi_{\an,\frac gf}\, (B)
\lr\frac{\ln\kn}{\kn}\rr^{1/2}\rb,
\end{equation*}
where $A$ and $B$ are some positive constants, which  depend only
on $c$.
\end{lemma}
\medskip
Lemma~\ref{lem_A.1} is a Bahadur-Kiefer type result. For a special case when $0<\alp<1$ is fixed it is stated in lemmas~3.1~\cite{gh2006} (cf.~also lemmas~4.1,~\cite{gh2007} and  Reiss~\cite{r}). We prove this lemma below in this section.

Note that we prove our results assuming that  $\alp_n<1$ for all sufficiently large $n$, where  $\alp_n$ can tend to $0$. Certainly, the same results can be obtained in case when $\alp_n>0$ for all sufficiently large $n$, where  $\alp_n$, in particular, can tend to $1$, i.e. on the right tail of the sample. Some new results on the Bahadur -- Kiefer representations for intermediate sample quantiles can be found in our recent paper \cite{gh2011}.

Lemma~\ref{lem_A.2} extends  lemma 4.3 from \cite{gh2007} (cf. also lemma~3.2,~\cite{gh2006}), where it was proved for
a~fixed $\alp$ to the~case that  $\an$ is a~sequence.
\begin{lemma}
\label{lem_A.2}
 Suppose that the conditions of
Lemma \ref{lem_A.1} hold true. Then
\begin{equation}
\label{(A.2)}
\int_{\qan}^{\qa}(G(x)-G(\qa))\, d\, F_n(x)=-\frac
12[F_n(\qa)-F(\qa)]^2\frac gf(\qa)+R_n,
\end{equation}
where  $\boldsymbol{P}(|R_n|>\Delta_n)=O\lr\kn^{-c}\rr$ for each $c>0$, and
\begin{equation*}
\Delta_n=A\, \frac{\an \ln\kn}{n}\lb \frac
{|g|}f(\qa)\lr\frac{\ln\kn}{\kn}\rr^{1/4}+\Psi_{\an,\frac gf}\, (B)
\rb,
\end{equation*}
where $A$,~$B$ are some positive constants, which  depend only
on $c$.
\end{lemma}
\begin{remark}{\sl
Suppose that $\kn^{-1}\ln n \to 0$, as $\nty$, and replace $\ln\kn$ by $\ln n$ in definition of function $\Psi_{\an,h}\, (B)$ (cf.~\eqref{1.8}) . Then lemmas~\ref{lem_A.1}, \ref{lem_A.2} remain  valid if we replace  $\ln\kn$ by $\ln n$ in formula for $\Delta_n$ in \eqref{(A.1)}--\eqref{(A.2)}. Furthermore, $\boldsymbol{P}(|R_n|>\Delta_n)=O\lr n^{-c}\rr$ for each $c>0$ in \eqref{(A.1)}--\eqref{(A.2)}.  To see the validity of this remark, it is enough to replace $\ln\kn$ by $\ln n$ in the proof of lemmas~\ref{lem_A.1} and ~\ref{lem_A.2} and use the assumption $\kn^{-1}\ln n \to 0$, no more changes in the proofs are needed. This remark is useful for obtaining of some results similar to Theorems~1.1, 1.4 and 1.6 in the case of light tails ($F$ has a~finite variance), it allows us to get the bounds of the order $O(n^{-r})$,  $0<r\le 1/2$, which are as one would expect in this case.}
\end{remark}
Let $U_1,\dots,U_n$ denote a sample of independent uniform $(0,1)$
distributed r.v.'s, and $U_{1:n}\le\cdots\le U_{n:n}$ -- the
corresponding order statistics. Put
\begin{equation}
\label{(A.3)}
N_{\an}^x = \sharp \{i : X_i \le \xi_{\an} \}\, ,\quad N_{\an} = \sharp
\{i : U_i \le \an \},
\end{equation}
and note that $\qan =X_{\kn:n}$ (because $\an=\kn/n$).

\medskip
\noindent{\bf Proof of lemma~\ref{lem_A.1}} We must prove that
$\boldsymbol{P}(|R_n|>\Delta_n)=O\lr\kn^{-c}\rr$ for each $c>0$ (cf. \eqref{(A.1)}),
and since the joint distribution of $X_{\kn :n}$, $N_{\an}^x $
coincide with joint distribution of $F^{-1}(U_{\kn :n})$,
$N_{\an}$ it is suffices to verify it for a remainder given by
\begin{equation*}
R_n= G(F^{-1}(U_{\kn :n}))-G(F^{-1}(\an))+\frac{N_{\an}-\an
n}{n}\frac gf(\qa).
\end{equation*}
Since $\boldsymbol{P}(U_{\kn :n}\notin U_a)=O(exp(-\delta n))$ for some
$\delta >0$ not depending on $n$, we can rewrite $R_n$ for all
sufficiently large $n$ as
\begin{equation}
\label{(A.4)}
\frac gf(\qa)\, R_{n,1}+ R_{n,2},
\end{equation}
where $R_{n,1}=U_{\kn :n}-\an+\frac{N_{\an}-\an n}{n}$, and
$R_{n,2}=\Bigl( \frac gf\lr F^{-1}(\an+\theta(U_{\kn :n}-\an
))\rr$ $-\frac gf \lr F^{-1}(\an)\rr\Bigr)\bigl(U_{\kn
:n}-\an\bigr)$, $0<\theta <1$. Fix an arbitrary $c>0$ and note
that we can estimate $R_{n,j}$, $j=1,2$, on the set $E=\bigl\{
\omega : |N_{\an}-\an n|<A_0\bigl(\an \, n\, \ln\kn \bigr)^{1/2}
\bigr\}$, where $A_0$ is a positive constant, depending only on
$c$, because by Bernstein inequality $\boldsymbol{P}(\Omega \setminus
E)=O(\kn^{-c})$ (in fact we can take every $A_0$:  $A^2_0>2c$). We will prove that
\begin{equation}
\label{(A.5)}
\boldsymbol{P}\bigl( |R_{n,1}|>A_1(\, \an\, )^{1/4}(\ln\kn/n)^{3/4}\bigr)
=O(\kn^{-c})
\end{equation}
and that
\begin{equation}
\label{(A.6)}
\boldsymbol{P}\bigl( |R_{n,2}|>A_2\, \an\, \Psi_{\an,\frac
gf}(B)(\ln\kn/\kn)^{1/2}\bigr) =O(\kn^{-c}).
\end{equation}
Here and elsewhere $A_i\, , i=1,2,\dots $, and $B$ denote some positive
constants, depending only on $c$. Relations \eqref{(A.4)}--\eqref{(A.6)} imply \eqref{(A.1)}.

First we prove \eqref{(A.5)} using a similar conditioning on $N_{\an}$
argument as in proof of lemmas 4.1, 4.3 in~\cite{gh2007}. First let $\kn\le N_{\an}$, then conditionally on
$N_{\an}$ the order statistic $U_{\kn:n}$ is distributed as
$\kn$-th order statistic $U'_{\kn:N_{\an}}$ of the sample
$U'_1,\dots,U'_{N_{\an}}$ independent $(0,\an)$ uniformly
distributed r.v.'s. Its expectation
$\boldsymbol{E}\big(U_{\kn:n}\lv\rp N_{\an},\ \kn\le
N_{\an}\big)=\an\, \frac{\kn}{N_{\an}+1}$, and the conditional
variance $V^2_{\kn}=\frac{\alp_n^2}{N_{\an}+2}\,
\frac{\kn}{N_{\an}+1}\, \bigl( 1-\frac{\kn}{N_{\an}+1}\bigr)$, and
on the set $E$ we have an estimate $V^2_{\kn}\le A_0(\an)^{1/2}\,
n^{-3/2}\ln^{1/2}\kn$. \ \ Then rewrite $R_{n,1}$ (at the event $\kn\le N_{\an}$) as
\begin{equation}
\label{(A.7)}
U_{\kn:n}-\an\, \frac{\kn}{N_{\an}+1}+R'_{n,1},
\end{equation}
where $R'_{n,1}=\an\, \frac{\kn}{N_{\an}+1}-\an+\frac{N_{\an}-\an
n}{n}=\frac{(N_{\an}-\kn)^2}{n(N_{\an}+1)}+\frac{N_{\an}-\kn}{n(N_{\an}+1)}
-\frac{\kn}{n(N_{\an}+1)}$, and on the set $E$ the latter
quantity is of the order $O\lr\frac{\ln\kn}{n}\rr$, and since
$\frac{\ln\kn}{n}=o\bigl( (\an)^{1/4}\lr
\frac{\ln\kn}{n}\rr^{3/4}\bigr)$, the remainder term $R'_{n,1}$
is of negligible order for our purposes. For the first two terms
in \eqref{(A.7)} we have
\begin{equation*}
\boldsymbol{P}\lr \lv U_{\kn:n}-\an\, \frac{\kn}{N_{\an}+1}\rv>A_1
(\an)^{1/4}\lr \frac{\ln\kn}{n}\rr^{3/4}\big| N_{\an}:\ \kn\le
N_{\an}\rr
\end{equation*}
\begin{equation}
\label{(A.8)}
=\boldsymbol{P}\lr \lv U'_{\kn:N_{\an}}-\an\, \frac{\kn}{N_{\an}+1}\rv>A_1
(\an)^{1/4}\lr \frac{\ln\kn}{n}\rr^{3/4}\rr =P_1+P_2, \quad
\end{equation}
where $N_{\an}$ is fixed, $\kn\le
N_{\an}$, $A_1$ is a constant which we will choose later,
$P_1=\boldsymbol{P}\lr U'_{\kn:N_{\an}}> \an\, \frac{\kn}{N_{\an}+1}+A_1
(\an)^{1/4}\lr \frac{\ln\kn}{n}\rr^{3/4}\rr$ and $P_2=\boldsymbol{P}\lr U'_{\kn:N_{\an}}< \an\, \frac{\kn}{N_{\an}+1}-A_1
(\an)^{1/4}\lr \frac{\ln\kn}{n}\rr^{3/4}\rr$.
We evaluate $P_1$, the treatment for $P_2$ is similar. Consider a
binomial r.v.
$S'_n=\sum_{i=1}^{N_{\an}}\bf{1}_{\{U'_{\, i:N_{\an}}\le \an\,
\frac{\kn}{N_{\an}+1}+A_1 (\an)^{1/4}\lr
\frac{\ln\kn}{n}\rr^{3/4} \}}$ with parameter $(p'_n,N_{\an})$,
where $p'_n=\min\bigl(1,\frac{\kn}{N_{\an}+1}+t_n \bigr)$, where $t_n=A_1 \lr
\frac{\log k_n}{\kn}\rr^{3/4} $. If $p'_n=1$, then $P_1=0$
and the inequality we need is valid trivial. Let $p'_n<1$ and let $\overline{S'}_n$ denote the average $S'_n/N_{\an}$, then  the probability  $P_1$ is equal to
\begin{equation}
 \label{(A.9)}
\boldsymbol{P}(S'_n<\kn)= \boldsymbol{P}\lr\overline{S'}_n-p'_n<\frac{k_n}{N_{\an}}-\frac{k_n}{N_{\an}+1} -t_n
\rr.
\end{equation}
Note that $\frac{k_n}{N_{\an}}-\frac{k_n}{N_{\an}+1}=\frac{k_n}{N_{\an}(N_{\an}+1)}<\frac 1{N_{\an}}$, and since  the latter quantity is $o\lr t_nk_n^{-1/4}\rr=o(t_n)$ on the set $E$,  this term  can be omitted at the r.h.s. of
\eqref{(A.9)} in our estimating. To evaluate  $ \boldsymbol{P}\lr\overline{S'}_n-p'_n< -t_n\rr$
we note that $p'_n-t_n=\frac{k_n}{N_{\an}+1}\in (0,1)$, and that $p'_n>1/2$ for all sufficiently large $n$ (and hence $k_n$ and  $N_{\an}$) on the set $E$. So, we may apply an inequality (2.2) of Hoeffding \cite{ho} with $\mu=p'_n$ and with $g(\mu)=1/(2\mu(1-\mu))$. Then we obtain
\begin{equation}
\label{(A.10)}
\boldsymbol{P}(S'_n<\kn)\le \exp \lr-N_{\an}t_n^2g(p'_n)\rr = \exp \lr-\frac{N_{\an}A_1^2 \bigl(\log k_n/\kn\bigr)^{3/2}}{2p'_n(1-p'_n)}\rr.
\end{equation}
Finally we note that  $1-p'_n=1-\frac{\kn}{N_{\an}+1}-A_1 \lr
\frac{\log k_n}{\kn}\rr^{3/4}\le \frac{N_{\an}+1-\kn}{N_{\an}+1}$, and on the set $E$ the latter quantity is not greater than $\frac{A_0(\kn\log k_n)^{1/2}}{N_{\an}}$. Then we can get a low bound for the ratio at the r.h.s. in \eqref{(A.10)}: $\frac{N_{\an}A_1^2 \bigl(\log k_n/\kn\bigr)^{3/2}}{2p'_n(1-p'_n)} \ge \frac{A_1^2 N^2_{\an} \bigl(\log k_n/\kn\bigr)^{3/2}}{2A_0(\kn\log k_n)^{1/2}}=\frac{A_1^2}{2A_0} \log k_n\lr\frac{N_{\an}}{k_n}\rr^2=\frac{A_1^2}{2A_0} \log k_n\lr 1+o(1)\rr$. This bound and \eqref{(A.10)} together
yield that when $\frac{A_1^2}{2A_0}\ge c$ the desired  relation  $P_1=O(k_n^{-c})$ hold true. The same estimate is valid for $P_2$.

In case $N_{\an}<\kn$ we use the fact that $U_{\kn:n}$
conditionally on $N_{\an}$ is distributed as $(\kn-N_{\an})$-th
order statistic $U''_{\kn-N_{\an}:n-N_{\an}}$ of the sample
$U''_1,\dots,U''_{n-N_{\an}}$ from $(1-\an,1)$ uniform
distribution, its expectation is
$\an+\frac{\kn-N_{\an}}{n-N_{\an}+1}$, and for the conditional
variance we have the estimate $V^2_{\kn-N_{\an}}\le
A_0(\ln\kn\an)^{1/2}\, n^{-3/2}$. In this case we use a
representation for $R_{n,1}=R''_{n,1}+R''_{n,2}$, where
$R''_{n,1}=U_{\kn:n}-\an-\frac{\kn-N_{\an}}{n-N_{\an}+1}(1-\an)$,
and $R''_{n,2}=\frac{N_{\an}-\an\,
n}{n}+\frac{\kn-N_{\an}}{n-N_{\an}+1}(1-\an)$. Similarly as in
first case we obtain that
$R''_{n,2}=O\bigl(\frac{\ln\kn}{n}\bigr)$ with probability
$1-O(\kn^{-c})$, and this term is of the negligible order in our
estimating. Using Hoeffding's inequality we obtain for $R''_{n,1}$
same estimate as for $R'_{n,1}$. So \eqref{(A.5)} is proved.

It remains to prove \eqref{(A.6)}. First note that by \eqref{(A.5)} on the set $E$
with probability $1-O(\kn^{-c})$ we have $|U_{\kn:n}-\an|\le
A_0\frac{(\kn\ln\kn)^{1/2}}{n}+A_1\an\lr\frac{\ln\kn}{\kn}\rr^{3/4}
=\bigl(\an\frac{\ln\kn}{n}\bigr)^{1/2}\bigl( 1+o(1)\bigr)$. Thus,
there exists $A_2$, depending only on $c$, such that $|R_{n,2}|\le
A_2 \bigl(\an\frac{\ln\kn}{n}\bigr)^{1/2}\Psi_{\an,\frac
gf}(A_2)$ with probability $1-O(\kn^{-c})$. This implies \eqref{(A.6)}.
The lemma is proved. $\quad \square$

\medskip
\noindent{\bf Proof of lemma~\ref{lem_A.2}} Let $N^x_{\an}$ and $N_{\an}$
are given as in \eqref{(A.3)}, then we can rewrite integral on the l.h.s. of
\eqref{(A.2)} as $\frac {sgn(N^x_{\an}-\kn )}{n}\sum_{i=(\kn\wedge N^x_{\an})+1}^{\kn\vee N^x_{\an}}(G(X_{i:n})-G(\qa))$, where $sgn(x)=x/|x|$, $sgn(0)=0$. Let us adopt the following notation: for
any integer $k$ and $m$ define a~set $I_{(k,m)}:=\{i: (k\wedge m)+1\le i\le k\vee m\}$ and let
$\sum_{i\in I_{(k,m)}}(.)_i:=sgn(m-k)\sum_{i=(k\wedge m)+1}^{k\vee m}(.)_i$. Then we must estimate
$R_n=\frac
1n\sum_{i\in I_{(\kn,{N^x_{\an})}}}(G(X_{i:n})-G(\qa))+\frac{\bigl(N^x_{\an}-\an n\bigr)^2}{2n^2}\frac gf(\qa)$ (cf. \eqref{(A.2)}), and
similarly as in proof of lemma \ref{lem_A.1} we note that $R_n$ is
distributed as
\begin{equation*}
\frac 1n\sum_{i\in I_{(\kn,{N_{\an})}}}\Bigl(G\circ F^{-1}(U_{i:n})-G\circ
F^{-1}(\an)\Bigr)+\frac{\bigl(N_{\an}-\an n\bigr)^2}{2n^2}\frac
gf(\qa) \phantom{R_{n,1}+R_{n,2}}
\end{equation*}
\begin{equation}
\label{(A.11)}
\phantom{R_{n,1}+R_{n,2}R_{n,1}+R_{n,2}+R_{n,1}}=\frac gf(\qa)R_{n,1}+R_{n,2},
\end{equation}
where $R_{n,1}=\frac 1n\sum_{i\in I_{(\kn,{N_{\an})}}}(U_{i:n}-\an)+\frac{\bigl(N_{\an}-\an n\bigr)^2}{2n^2}$,
$R_{n,2}=\frac 1n\sum_{i\in I_{(\kn,{N_{\an})}}}\lb \frac gf\circ
F^{-1}\bigl(\an +\theta_i(U_{i:n}-\an)\bigr)-\frac gf\circ
F^{-1}\bigl(\an\bigr)\rb\Bigl(U_{i:n}-\an\Bigr)$, where
$0<\theta_i<1$, \ $i\in I_{(\kn,{N_{\an})}}$.

 Fix an arbitrary $c>0$ and prove that
\begin{equation}
\label{(A.12)}
\boldsymbol{P}\Bigl( |R_{n,1}|>A_1(\, \an\, )^{3/4}(\ln\kn/n)^{5/4}\Bigr)
=O(\kn^{-c}),
\end{equation}
\begin{equation}
 \label{(A.13)}
\boldsymbol{P}\Bigl( |R_{n,2}|>A_2\, \an\, \frac{\ln\kn}{n}\Psi_{\an,\frac
gf}(A_2)\Bigr) =O(\kn^{-c}).
\end{equation}
Since $\bigl(\, \an\,
\bigr)^{3/4}\bigl(\ln\kn/n\bigr)^{5/4}=\an\frac{\ln\kn}{n}
\bigl(\ln\kn/\kn\bigr)^{1/4}$, \ relations \eqref{(A.12)}--\eqref{(A.13)} imply
\eqref{(A.2)}. Note that as in proof of lemma \ref{lem_A.1} it is enough to estimate
$R_{n,j}$, $j=1,2$, on the set $E=\bigl\{ \omega : |N_{\an}-\an
n|<A_0\bigl(\an \, n\, \ln\kn \bigr)^{1/2} \bigr\}$, where
$A_0>0$ is a constant, depending only on $c$, such that $\boldsymbol{P}(\Omega
\setminus E)=O(\kn^{-c})$.

First we treat $R_{n,2}$. Note that
\begin{equation*}
\max_{i\in I_{(\kn,{N_{\an})}}}\big|\, U_{i:n}-\an\big|=\big|\,
U_{\kn:n}-\an\big|\vee \big|\, U_{N_{\an}:n}-\an\big|\vee \big|\, U_{N_{\an}+1:n}-\an\big|\, ,
\end{equation*}
$\boldsymbol{P}\Bigl(\big|\, U_{\kn:n}-\an\big|>A_0\bigl(\an\ln\kn /n\,
\bigr)^{1/2}\Bigr)=O(\kn^{-c})$ (cf. proof of lemma~\ref{lem_A.1}), and for
$j=N_{\an:n}\, , N_{\an:n}+1$ simultaneously we have  $\boldsymbol{P}\Bigl(\big|\,
U_{j:n}-\an\big|>A_1\frac{\ln\kn}{n}\Bigr)\le \boldsymbol{P}\Bigl(
U_{N_{\an}+1:n}-U_{N_{\an}:n}>A_1\frac{\ln\kn}{n}\Bigr) =\boldsymbol{P}\Bigl(U_{1:n}>A_1\frac{\ln\kn}{n}\Bigr)=
\Bigl(1-A_1\frac{\ln\kn}{n}\Bigr)^n=O(\kn^{-c})$ for $A_1>c$.
Since $\frac{\ln\kn}{n}=o(\frac{\an\ln\kn}{n})^{1/2}$, \ on the
set $E$ we obtain
\begin{equation*}
\big|R_{n,2}\big|\le \frac 1n\Psi_{\an\frac
gf}(A_0)A_0^2\Bigl(\an\, n\ln\kn\Bigr)^{1/2}\Bigl(\frac{\an\,
\ln\kn}{ n}\Bigr)^{1/2}=A_2\an\frac{\ln\kn}{n}\Psi_{\an\frac
gf}(A_0)
\end{equation*}
with probability $1-O(\kn^{-c})$, and \eqref{(A.13)} is proved.

Finally, consider $R_{n,1}$. Note that conditionally on $N_{\an}$,
\ $\kn\le N_{\an}$, the order statistics $U_{i:n}$, $\kn\le i\le N_{\an}$, are distributed
as the order statistics $U'_{i:N_{\an}}$
from the uniform $(0,\an)$ distribution (cf.~lemma~5.1, Section~5), their conditional
expectations are equal to $\an\frac{i}{N_{\an}+1}$. Then in the case $\kn\le  N_{\an}$ (the proof for the~case $N_{\an}<\kn$ is similar
(cf. proof of lemma~\ref{lem_A.1}) with respect to interval $(1-\an,1)$, and
we omit the details) we rewrite $R_{n,1}$ as
\begin{equation}
\label{(A.14)}
R_{n,1}=\frac
1n\sum_{i=\kn+1}^{N_{\an}}\bigl(U_{i:n}-\an\frac{i}{N_{\an}+1}\bigr)+R'_{n,1}\,
,
\end{equation}
where $R'_{n,1}=\frac
1n\sum_{i=\kn+1}^{N_{\an}}\an\bigl(\frac{i}{N_{\an}+1}-1\bigr)+ \frac{(N_{\an}-\an\,
n)^2}{2n^2}=-\frac{\kn}{n^2}\frac{(N_{\an}-\kn)(N_{\an}-\kn-1)}{2\,(N_{\an}+1)}+ \frac{(N_{\an}-\kn)^2}{2n^2}
=\frac{(N_{\an}-\kn)^2(N_{\an}+1-\kn)}{2\,(N_{\an}+1)n^2}-\frac{\kn(N_{\an}-\kn)}{2(N_{\an}+1)n^2}$, and on the set $E$
the latter quantity is of the order $O\lr
\frac{\kn^{1/2}(\ln\kn)^{3/2}}{n^2}\rr=o\lr\bigl(\an
\bigr)^{3/4}\lr\frac{\ln\kn}{n}\rr^{5/4}\rr$,  i.e. $R'_{n,1}$ is
of negligible order (cf.~\eqref{(A.12)}) for our purposes.

It remains to evaluate the dominant first term on the r.h.s. in
\eqref{(A.14)}. Fix an arbitrary $c_1>c+1/2$, and note that conditional on
$N_{\an}$ the variance of $U_{i:n}$ ($\kn+1\le i\le N_{\an}$) is equal
to $V_i^2=\bigl(\an\bigr)^2\frac 1{N_{\an}+2}\,
\frac{i}{N_{\an}+1}\, \lr 1-\frac{i}{N_{\an}+1}\rr$, and on the
set $E$ it is less than
$\bigl(\an\bigr)^2\frac{A_0\kn^{1/2}(\ln\kn)^{1/2}}{N^2_{\an}}$,
and $V_i\le\an A_0^{1/2}\kn^{1/4}(\ln\kn)^{1/4}/N_{\an}
\le A_0^{1/2}\an \kn^{-3/4}(\ln\kn)^{1/4}\le A_0^{1/2}\bigl(\an\bigr)^{1/4}n^{-3/4}(\ln\kn)^{1/4}$. Using
Hoeffding's inequality (similarly as in proof of lemma~\ref{lem_A.1}), we find
that
\begin{equation*}
\boldsymbol{P}\Bigl( \big|U_{i:n}-\an\frac{i}{N_{\an}+1} \big|>A_1
\bigl(\an\bigr)^{1/4}\bigl(\ln\kn/n \bigr)^{3/4}\Big| N_{\an}:\
\kn\le N_{\an}\Bigr)=O(\kn^{-c})\, ,
\end{equation*}
where $A_1$ depends only on $c_1$ (in fact it is true for every
$A_1$ such that $A^2_1>2A_0c$). Thus,
\begin{equation*}
\boldsymbol{P}\Bigl( \frac 1n
\big|\sum_{i=\kn}^{N_{\an}}\bigl(U_{i:n}-\an\frac{i}{N_{\an}+1}
\bigr)\big|>A_0A_1 \bigl(\an\bigr)^{3/4}\bigl(\ln\kn/n
\bigr)^{5/4}\Big| N_{\an}:\ \kn\le N_{\an}\Bigr)
\end{equation*}
\begin{equation}
\label{(A.15)}
\qquad \qquad \qquad \qquad \qquad
\qquad\qquad \quad \ \le A_0(\kn\ln\kn )^{1/2}\,
O(\kn^{-c_1})=O(\kn^{-c})\, .
\end{equation}
Combining \eqref{(A.14)}--\eqref{(A.15)} and similar estimates for the case
$N_{\an}<\kn$, we come to \eqref{(A.12)}. The lemma is proved.$\quad \square$
\section{Appendix}
\label{appendix}
Let as before, $\na=\sharp \{i \ : X_i \le \xia ,\ i=1,\dots,n\}$, where $0<\alp <1$ is fixed. In this appendix we prove that
conditionally on $\na$ the order statistics $X_{1:n},\dots ,X_{\na:n}$ are distributed as
order statistics corresponding a sample of $\na$ i.i.d.~r.v.'s
with distribution function $F(x)/\alp$, $x\le \xia$. Though this fact is known (cf.~\cite{kallen},~\cite{s}), we give a~brief proof of it. Let
$U_1,\dots,U_n$ be independent r.v.'s uniformly distributed on
$(0,1)$ and let $U_{1:n},\dots,U_{n:n}$ denote the corresponding
order statistics. Put $\nau=\sharp \{i \ : \ U_i \le \alp \}$.
Since $\xin \stackrel{\rm d}{=} F^{-1}(\uin)$ and $\na
\stackrel{\rm d}{=} \nau\,$, it is enough to prove the assertion
for the uniform
distribution.
\begin{lemma}
\label{Lemma_A} Conditionally given $\nau$, the order
statistics $U_{1,n},\dots,U_{\nau,n}$ are distributed as order
statistics corresponding to a sample of $\nau$ independent
$(0,\alp)$-uniform distributed r.v.'s.
\end{lemma}
\noindent{\bf Proof}.  $a).$ First consider the case $\nau= n$. Take
arbitrary $0<u_1\le \cdots \le u_n < \alp$ and write
\begin{equation*}
\label{ap_1}
\begin{split}
P(U_{1:n} \le u_1,\dots,U_{\nau:n}\le u_n \mid\nau =n)
=\frac{P(U_{1:n} \le u_1,\dots,U_{n:n}\le u_n)}{\alp^n}\\
=\frac{n!}{\alp^n}\int_0^{u_1}\int_{u_1}^{u_2}\dots
\int_{u_{n-1}}^{u_n} d\,x_1 d\,x_2\dots d\,x_n,
\end{split}
\end{equation*}
and the latter is $d.f.$ of the order statistics corresponding to
the sample of \ $n$ \ independent $(0,\alp)$-uniform distributed
r.v.'s. \ \  $b).$ Consider the case $\nau=k<n$. Let $F_{i,n}(u)=P(U_{i:n}\le u)$ be a~$df$ of $i$-th order statistic, put $P_n(k)=P(\nau=k)={n\choose k} \alp^k(1-\alp)^{n-k}$. Then we have
\begin{equation}
\label{ap_2}
P(U_{1:n} \le u_1,\dots,U_{\nau :n}\le u_k \mid \nau =k)
= \frac {P(U_{1:n} \le u_1,\dots,U_{k :n}\le u_k, U_{k+1:n}>\alp)}{P_n(k)}.
\end{equation}
The probability in the nominator on the r.h.s. of  \eqref{ap_2} is equal to
\begin{equation*}
\int_{\alp}^1 \!\!\! P\big(U_{1:n} \le
u_1,\dots,U_{k:n} \le u_k \mid U_{k+1:n}=v\big)\,d F_{k+1,n}(v),
\end{equation*}
and by the Markov property of order statistics the latter quantity equals
\begin{equation*}
\label{ap_3}
\begin{split}
&\phantom{ =}\int_{\alp}^1 \lr \frac{k!}{v^k} \int_0^{u_1}\int_{u_1}^{u_2}\dots \int_{u_{k-1}}^{u_k} d\,x_1
d\,x_2\dots d\,x_k\rr \,d F_{k+1,n}(v)\\
&= \frac{k!}{\alp^k}\lr \int_0^{u_1}\int_{u_1}^{u_2}\dots \int_{u_{k-1}}^{u_k} d\,x_1
d\,x_2\dots d\,x_k \rr\times \alp^k\int_{\alp}^1 \frac 1{v^k} \,d F_{k+1,n}(v),
\end{split}
\end{equation*}
and since $\alp^k\int_{\alp}^1 \frac 1{v^k} \,d F_{k+1,n}(v)=\alp^k\int_{\alp}^1 \frac{(1-v)^{n-k-1}}{B(k+1,n-k)} \,d v
= {n\choose k}\alp^k (1-\alp)^{n-k}=P_n(k)$, where $B(k+1,n-k)=k!(n-k-1)!/n!$, we obtain that conditional probability in  \eqref{ap_2} is equal
\begin{equation*}
\frac{k!}{\alp^k}\int_0^{u_1}\int_{u_1}^{u_2}\dots
\int_{u_{k-1}}^{u_k} d\,x_1 d\,x_2\dots d\,x_k,
\end{equation*}
which  corresponds to the $(0,\alp)$-uniform distribution. The lemma
is proved.$\quad \square$


\begin{thebibliography} {99}
\bibitem{b}
Bahadur,~R.R.(1966). A note on quantiles in large samples., {\it
 Ann. Math. Statist.} {\bf 37}~577--580.
\bibitem{bgzit}
Bentkus,~V., G\"otze,~F. and  Zitikis R. (1994). Lower
estimates of the convergence rate for $U$-statistics.,  {\it Ann.
Probab.} {\bf~22}~1707--1714.
\bibitem{bjz}
Bentkus,~V., Jing,~B-Y. and Zhou,~W.( 2009). On Normal Approximations to U -statistics., {\it Ann.
Probab.} {\bf~37}~2174--2199.
\bibitem{bgz}
Bentkus,~V., G\"otze,~F. and  van Zwet,~W.R.(1997). An
Edgeworth expansion for symmetric statistics., {\it Ann. Statist.} {\bf~25}~851--896.
\bibitem{bi}
Bickel,~P.J.,G\"otze,~F. and  van Zwet,~W.R. (1986).  The
Edgeworth expansion for $U-$ statistics of degree two., {\it Ann.
Statist.} {\bf~14}~1463--1484.
\bibitem{bingham}
Bingham,~N.M., Goldie,~C.M. and Teugels,~J.L. (1987)., {\it Regular variation}, Cambridge Univ. Press (Encyclopedia
Math. Appl.) {\bf 27}, Cambridge.
\bibitem{bormogu}
Borovkov,~A.A. and  Mogulskii,~A.A. (2006). On large and
superlarge deviations of sums of independent random vectors under
Cramer's condition.~II., {\it Theory Probab. Appl.} {\bf~51}~641--673.
\bibitem{cs}
Chen,~L.H.Y. and  Shao,~Q.M. (2007) Normal approximations for nonlinear statistics using a concentration inequality approach., {\it Bernoulli} {\bf~13}~581--599.
\bibitem{chm}
Cs\"org\H{o},~S., Haeusler,~E. and  Mason,~D.M. (1988) The
asymptotic distribution of trimmed sums., {\it Ann. Probab.} {\bf~16}~672--699.
\bibitem{cm}
Cs\"org\H{o},~S. and  Megyesi,~Z. (2001) Trimmed sums from the domain of geometric partial attraction of semistable laws. In:  State of the art in probability and statistics: Festschrift for Willem R. van Zwet (M.de Gunst  C.Klaassen  A.van der Vaart,~eds), {\it  Lecture Notes - Monograph Series} {\bf~36}~173--194, Institute of Mathematical Statistics, Beachwood, Ohio.
\bibitem{en}
Egorov,~V.A. and  Nevzorov,~V.B. (1974) Certain estimates of the
rate of convergence of sums of order statistics to the normal law., {\it Zap. Nauch. Sem. LOMI/ Leningrad. otdel. math. inst. V.A.Steklov} {\bf~41}~105--128 (in Russian); (1978). Transl. in {\it J. Math. Sci. (New York)} {\bf~9}, No.~1,~81--105.
\bibitem{feller}
Feller,~W. (1971). {\it An introduction to probability
theory and its applications}, vol. II, John Wiley \& Sons, New
York.
\bibitem{friedrich}
Friedrich,~K.O. (1989).  A Berry -- Esseen bound for
functions of independent random variables., {\it Ann. Statist.} {\bf~17}~170--183.
\bibitem{gri}
Gribkova,~N.V. (1993).  On analogs of the Berry -- Esseen
inequality for truncated linear combinations of order statistics., {\it  Theory Probab. Appl.} {\bf~38}~142--149.
\bibitem{gh2006}
Gribkova,~N.V. and Helmers,~R. (2006). The empirical
Edgeworth expansion for a Studentized trimmed mean., {\it Math.~Methods~Statist.} {\bf~15}~61--87.
\bibitem{gh2007}
Gribkova,~N.V. and Helmers,~R. (2007). On the Edgeworth
expansion and the $M$ out of $N$ bootstrap accuracy for a
Studentized trimmed mean.,{\it Math.~Methods~Statist.} {\bf~16}~142--176.
\bibitem{gh2007a}
Gribkova,~N.V. and Helmers,~R. (2010). On the consistency of the $M \ll N$ bootstrap approximation for a trimmed mean., {\it Theory Probab. Appl.} {\bf~55}, no.~1,~42--53.
\bibitem{gh2011}
Gribkova,~N.V. and Helmers,~R. (2011).  On the Bahadur -- Kiefer representation for intermediate sample quantiles {\it (submitted paper)}.
\bibitem{griffin}
Griffin,~P.S. and Pruitt,~W.E. (1989) Asymptotic normality
and subsequential limits of trimmed sums., {\it Ann. Probab.} {\bf~17}~1186--1219.
\bibitem{hp}
Hall,~P. and  Padmanabhan,~A.R. (1992) On the bootstrap and
the trimmed mean., {\it J. of Multivariate Analysis} {\bf~41}~132--153.
\bibitem{h}
Helmers,~R. (1991)  On the Edgeworth expansion and the
bootstrap approximation for a Studentized $U-$statistic., {\it Ann.
Statist.} {\bf~19}~470--484.
\bibitem{hjq}
Helmers,~R., Jing,~B.-Y.,  Qin,~G. and Zhou,~W. (2004) Saddlepoint approximations  to the trimmed mean., {\it Bernoulli} {\bf~10}, no.~3,~465--501.
\bibitem{ho}
Hoeffding,~W. (1963). Probabilities inequalities for sum
of bounded random variables., {\it Amer.Stat.Assoc.J.} {\bf 58} 13-30.
\bibitem{kallen}
Kallenberg,~O. (2002)., {\it Foundations of Modern Probability}, Springer,
New York.
\bibitem{kiefer}
Kiefer,~J.C. (1970) Deviations between the sample quantile
process and the sample $df$.  In: {\it Nonparametric Techniques in
Statistical Inference} (M. Puri, ed.)~299--319, London, Cambridge
Univ. Press.
\bibitem{kulik}
Kulik,~R. (2008). Trimmed sums of long range dependent moving averages., {\it Statistics~\&~Probability Letters} {\bf~78}~2536--2542.
\bibitem{peng}
Peng,~L. (2001). Estimating the mean of a heavy tailed
distribution., {\it Statistics~\&~Probability Letters} {\bf~52}~255--264.
\bibitem{petrov}
Petrov,~V.V. (1975). {\it  Sums of independent random
variables}, Springer-Verlag, New York.
\bibitem{pz}
Putter,~H. and  van Zwet,~W.R. (1998).  Empirical Edgeworth
expansions for symmetric statistics., {\it Ann. Statist.} {\bf~26}~1540--1569.
\bibitem{r}
Reiss,~R.-D. (1989). {\it Approximate distributions of order
statistics with applications to nonparametric statistics},
Springer-Verlag, New York.
\bibitem{shorack}
Shorack,~G.R. and  Wellner,~J.A. (1986). {\it Empirical
processes with application in statistics}, Wiley, New York.
\bibitem{s}
Stigler,~S.M. (1973).  The asymptotic distribution of the
trimmed mean., {\it Ann. Statist.} {\bf~1}~472--477.
\bibitem{zwet}
van Zwet,~W.R. (1984). A Berry -- Esseen bound for symmetric
statistics., {\it Z. Wahrsch. Verw. Gebiete} {\bf~66}~425--440.
\end{thebibliography}
\end{document}